\documentclass[letterpaper, 10 pt, conference]{ieeeconf}  

\IEEEoverridecommandlockouts                              
\usepackage{amsthm}

\def\QEDclosed{\mbox{\rule[0pt]{1.3ex}{1.3ex}}} 

\def\QED{\QEDclosed} 

\newtheoremstyle{new_plain}
  {0.5\topsep}   
  {0.5\topsep}   
  {\itshape}  
  {0em}       
  {\bfseries} 
  {.}         
  {5pt plus 1pt minus 1pt} 
  {}          
\theoremstyle{new_plain}
\newtheorem{theorem}{Theorem}
\newtheorem{lemma}{Lemma}
\newtheorem{assumption}{Assumption}
\newtheorem{corollary}{Corollary}

\newtheoremstyle{new_definition}
  {0.5\topsep}   
  {0.5\topsep}   
  {\normalfont}  
  {0pt}       
  {\bfseries} 
  {.}         
  {5pt plus 1pt minus 1pt} 
  {}          
\theoremstyle{new_definition}
\newtheorem{definition}{Definition}

\newtheoremstyle{new_remark}
  {0.5\topsep}   
  {0.5\topsep}   
  {\normalfont}  
  {1em}       
  {\itshape} 
  {.}         
  {5pt plus 1pt minus 1pt} 
  {}          
\theoremstyle{new_remark}
\newtheorem{remark}{Remark}

\usepackage[dvipsnames]{xcolor}
\usepackage{cancel}
\usepackage{graphicx} 
\usepackage{epsfig} 
\usepackage{mathptmx} 
\usepackage{amsmath} 
\usepackage{amssymb}  
\usepackage{xcolor}
\usepackage[cal=cm,scr=boondox]{mathalpha}

\usepackage{stmaryrd}

\usepackage{subcaption}

\title{\LARGE \bf
Reference Replay for Excitation Injection and Subsequent Exponential Convergence under Concurrent Learning Adaptive Control
}

\author{\parbox{3 in}{\centering Tom Kaufmann and Johann Reger
\thanks{The authors are with the Control Engineering Group, Technische Universität Ilmenau, P.O. Box $10$ $05$ $65$, D-$98684$ Ilmenau, Germany}
\thanks{Corresponding author: {\tt\footnotesize johann.reger@tu-ilmenau.de}}}
}

\begin{document}
\addtolength{\skip\footins}{-1mm}
\addtolength{\abovedisplayskip}{-0.65mm}
\addtolength{\belowdisplayskip}{-0.65mm}
\maketitle
\thispagestyle{empty}
\pagestyle{empty}

\begin{abstract}
Concurrent learning (CL) guarantees precise identification of unknown parameters---if data with finite excitation (FE) happens to be collected online. To improve the reliability of CL, this work establishes a mechanism allowing FE-occurrence to be verified a priori: We propose replaying a given training segment in the reference until data with finite excitation is~col-lected online. Provided the training segment fulfills a certain~off-line excitation test, FE then occurs in unknown, but finite~time. Moreover, we provide an upper bound on the time-to-FE. Even\linebreak though this bound involves a dependency on the unknown initial Lyapunov-function value, its structure shows that FE occurs in arbitrarily early replay cycles when increasing the control gain; which may motivate employing a high-gain strategy. The theo-\linebreak retical findings are verified on an academic simulation example.
\end{abstract}

\vspace{-3.75mm}\section{Introduction}\label{section:INTRODUCTION}\vspace{-1.25mm}
Classical adaptive control facilitates trajectory tracking when faced with structured uncertainty in the dynamical behavior of a plant by providing estimates for unknown parameters that enter the dynamics linearly through some known, possibly nonlinear regressor. While classical adaptive control is capable of driving the tracking error toward zero and thereby achieves the tracking objective, the estimation error generally can only be guaranteed to remain bounded~\cite{narendra_stable_1989}. This is due to the classical update law of the estimation, as it renders the time derivative of the Lyapunov certificate for its tracking capability independent of the estimation error. That is, the estimates of classical, Lyapunov-based adaptation are designed solely to counteract the immediate effect of the structured uncertainty. As a consequence, these estimates are useless for the purpose of identifying the unknown parameters---unless the regressor evaluated on the reference trajectory happens to fulfill a restrictive persistency-of-excitation (PE) condition; see \cite{loria_-persistency_2002},~\cite{loria_uniform_2003} for details. 

The desire to relax the restrictive requirement of the classical, Lyapunov-based approach for PE spurs the development of data-driven strategies in the field of adaptive control. Specifically, the approaches \cite{chowdhary_concurrent_2010,chowdhary_exponential_2014,pan_model_2016,kamalapurkar_concurrent_2017,parikh_integral_2019,shen_regulation-triggered_2026}, which fall under the umbrella of concurrent learning (CL), propose to augment the classical update law of the estimation with terms built from data that is recorded during closed-loop operation. The common intention of all aforementioned CL techniques is to remember already extracted information about the unknown parameters in order to enable learning during periods where the reference might lack excitation. Notably, the part of the time derivative of the Lyapunov function that originates from the augmentation of the update law with data is negative semi-definite w.r.t. the estimation error---with   the potential of becoming negative definite once data with finite excitation (FE) has been collected. As CL relies on the regressor evaluated on the closed-loop trajectory of the actual state, the occurrence of FE has to be verified online. However, in case the collected data fulfills FE, the time derivative of the Lyapunov function becomes negative definite w.r.t. the estimation error; which then guarantees~the precise identification of the unknown parameters---even if afterward, the reference trajectory lacks the excitation that would have been needed with the classical data-free approach to identify these parameters. Moreover, even for tasks where identification of the unknown parameters is not requested, collecting data remains desirable because CL with a quadratic Lyapunov function, subsequent to observing FE in closed loop, does guarantee the \textit{exponential} convergence toward zero of the estimation error as well as of the tracking error. Thereby, adaptation with CL under FE offers an~up-\linebreak grade over the guarantee for mere \textit{asymptotic} tracking that can be given under the classical data-free approach, as the proof for the tracking capability of the latter without~PE is reliant on Barbalat's lemma (or one of its multiple alternatives neatly listed in \cite{sun_gathering_2023}, e.g.\ \cite{tao_simple_1997}).

In this work, we mitigate the severity of one specific~draw-back of CL: Its appealing capability to render estimation and tracking errors exponentially convergent toward zero hinges on the occurrence of FE in closed loop. To~guarantee\linebreak closed-loop FE-occurrence a priori, we propose replaying some pre-defined, sufficiently exciting training segment in the reference trajectory until FE is injected into the data that CL is collecting during closed-loop operation. Therein, the\linebreak training segment has to pass some a-priori-verifiable, offline\linebreak excitation test and the instants at which CL acquires data are required to include replay-periodic copies of training data-acquisition instants so that the training excitation carries over to the data collected on the actual closed-loop trajectory.\footnote{In \cite{obuz_adaptive_2024}, the CL data-augmentation of the update law (\cite{obuz_adaptive_2024}, Eq. (15))~is built from the regressor evaluated on the reference (as opposed to the~regressor on the actual trajectory, as we use herein in unison with the established CL literature \cite{chowdhary_concurrent_2010}-\cite{shen_regulation-triggered_2026}), i.e. with FE supposedly verifiable a priori. They claim exponential convergence of the tracking and the estimation errors. However, their proof contains an unjustified bound: (\cite{obuz_adaptive_2024}, Eq.\ (11)) bounds each discrepancy term $\Omega(t_i)$ in terms of the tracking error $e(t_i)$ at the respective data-acquisition instants $t_i$. In (\cite{obuz_adaptive_2024}, Eq.\ (17)), this bound is converted into an estimate proportional to the current tracking error $e(t)$ without establishing a relation between historical $e(t_i)$ and current $e(t)$ errors. Consequently, the subsequent Lyapunov estimate (\cite{obuz_adaptive_2024}, Eq.\ (31)) and hence, the claimed exponential-convergence result (\cite{obuz_adaptive_2024}, Thm.\ 1) are, as of now, unproven.}
The asymptotic convergence of the tracking error toward zero then is instrumental for \textit{our central result} (see~Theorems\linebreak \ref{Theorem:ICL}, \ref{Theorem:NameWärNett} below) \textit{that, under suitable smoothness requirements, this injection of excitation happens in unknown, but provably finite time.} Transient modification of the reference to obtain FE is pursued also by Gallegos and Aguila-Camacho for an integral-memory estimator in $\!$\cite{gallegos_relaxed_2024}, where they temporarily~re-\linebreak place the reference by a signal on which the regressor is PE, until FE of the integral memory is detected online. Similarly,$\!$\footnote{A more detailed contextualization of reference replay w.r.t.\ \cite{gallegos_relaxed_2024} can be found in Subsection~\ref{subsection:RRvsTE}, where we discuss their similarities $\&$ differences.}\linebreak once FE of the collected data is observed during operation of our closed loop---as we can guarantee a priori to happen in finite time when using reference replay up until this time instant---any desired reference trajectory can be tracked afterwards. It is~nonetheless important to acknowledge that reference replay buys the appealing a-priori predictability of closed-loop FE-occurrence and hence, of precise parameter identification at the expense of postponing the optionality to track any freely chosen, albeit trackable,~desired~reference.

\textit{Notation.} In Section~\ref{section:RESULTS}, we use $\llbracket\cdot\rrbracket:[0,\infty)\to[0,\infty)$ to denote $\llbracket\mathscr{r}\rrbracket=\lfloor\mathscr{r}\rfloor$ for $\mathscr{r}\in[0,\infty)\setminus\mathbb{N}$ and $\llbracket\mathscr{r}\rrbracket=\mathscr{r}-1$ for $\mathscr{r}\in\mathbb{N}$. That is, $\llbracket\cdot\rrbracket$ acts on the nonnegative numbers like a floor-function that rounds down also the natural numbers to their integer predecessors. Otherwise, the notation is standard.

\vspace{-0.5mm}\section{Preliminaries}\label{section:PRELIMINARIES}\vspace{-0.5mm}
Before developing the reference-replay mechanism and then proving its properties  in Section~\ref{section:RESULTS}, we establish the necessary infrastructure and recall some preliminary results.

To focus this work on the specifics of the update law while also retaining design freedom due to a wide class of trackable references ($x_\mathrm{r}(t)\in\mathbb{R}^n$, $t\in[0,\infty)$ with~$x_\mathrm{r}\in C^1\cap L^\infty$, $\dot{x}_\mathrm{r}\in L^\infty$), we consider a plant with the simplistic~dynamical~behavior
\begin{align}
\label{Plant}
    \dot{x}&=Y(x)\theta+u.
\end{align}
The structure of~\eqref{Plant} is used also in \cite{parikh_integral_2019,gallegos_relaxed_2024, wen_relaxation_2026} in order to study different CL techniques and acts as a proxy that is~easi-ly adoptable to more involved model-reference structures~\cite{pan_model_2016}. The state therein, $x(t)\in\mathbb{R}^n$, is available for being fed back into the dynamical behavior via the input signal $u(t)\in\mathbb{R}^n$ with the aim of driving the tracking error $e=x-x_\mathrm{r}$ toward zero. Further, $Y(x)\theta$ represents the structured uncertainty in~\eqref{Plant}. It consists of the known, continuous regressor $Y:\mathbb{R}^n\to\mathbb{R}^{n\times p}$ and the unknown, but constant parameters $\theta\in\mathbb{R}^p$. Due to the unavailability of $\theta$, the adaptive control law
\begin{align}
\label{AdaptiveControl}
    u&=-Ke-Y(x)\hat{\theta}+\dot{x}_\mathrm{r},
\end{align}
where $0\prec K\in\mathrm{sym}(n)$ is some positive definite gain matrix, uses the estimates $\hat{\theta}(t)\in\mathbb{R}^{p}$ instead. Let $\tilde{\theta}=\hat{\theta}-\theta$ be the estimation error. This leads to the tracking-error dynamics~as
\begin{align}
\label{ErrorDyn}
    \dot{e}&=-Ke-Y(x)\tilde{\theta}.
\end{align}
The pending design of the update law $\dot{\hat{\theta}}$ of the estimation (which is restricted to available information) offers the opportunity to influence the evolution of the Lyapunov~function
\begin{align}
\label{LyapFcn}
    V(e,\tilde{\theta})&=\tfrac{1}{2}\|e\|^2+\tfrac{1}{2}\tilde{\theta}^\top\Gamma^{-1}\tilde{\theta}
\end{align} by purposefully manipulating its time derivative
\begin{align}
\label{dotV}
    \dot{V}(e,\tilde{\theta})&=-e^\top K e+\tilde{\theta}^\top\big(\Gamma^{-1}\dot{\hat{\theta}}-Y(x)^\top e\big).
\end{align}
Therein, $0\prec \Gamma\in\mathrm{sym}(p)$ is a freely chosen adaptation rate~of the different variations of the update law that we explore subsequently in the data-free and in the~data-augmented~setting.

\vspace{-.75mm}\subsection{Classical data-free update~law}\vspace{-.75mm}
The update law of the estimation is classically chosen as
\begin{align}
\label{UpdateLaw}
    \dot{\hat{\theta}}&=\Gamma\,Y(x)^\top e,
\end{align}
so that it cancels out the indefinite (and hence, parasitic) addend on the right-hand-side of \eqref{dotV}; thus rendering the time derivative of the Lyapunov function into $\dot{V}(\cdot)=-e^\top Ke$, i.e.\ negative definite w.r.t.\ the tracking error and independent of the estimation error. By means of employing standard arguments, see e.g. \cite{narendra_stable_1989}, one readily obtains {the following properties} of the closed loop that represents the plant \eqref{Plant} under adaptive control \eqref{AdaptiveControl} with the classical update law~\eqref{UpdateLaw}:\begin{enumerate}
    \item\label{Stability} $V(\cdot)$ is positive definite w.r.t.\ the full state $z=(e^\top,\tilde{\theta}^\top)^\top$ of the closed loop whereas its time derivative is nonpositive. Thereby, {the origin of $z$ is stable}.\allowdisplaybreaks
    \item\label{Boundedness} As $V(\cdot)$ is radially unbounded in $z$, its sublevel~set $\Omega_c=\{z\in\mathbb{R}^{n+p}:V(\cdot)\leq c\}$, at any given level $c>0$,~is bounded. Due to $\dot{V}(\cdot)\leq 0$, one has $V(\cdot)|_{t}\leq V(\cdot)|_{t=0}$ and therewith obtains $\forall t\geq 0: z(t)\in\Omega_{V(\cdot)|_{t=0}}$. That is, for any given (unknown) initial condition $z(0)$, $z$ results as a bounded function of time.~In particular, {this guarantees the estimates $\hat{\theta}=\tilde{\theta}+\theta$ to remain~bounded}.\allowdisplaybreaks
    \item\label{AsymptoticTracking} Given the fact $V(\cdot)|_{t}\!-\!V(\cdot)|_{t=0}=\int_{0}^t\!\dot{V}(\cdot)|_{s}\mathrm{d}s$~in~view of $V(\cdot)|_{t}\geq 0$, $\dot{V}(\cdot)\leq -\lambda_{\min}(K)\|e\|^2$ leads to~$V(\cdot)|_{t=0}\geq\allowbreak \lambda_{\min}(K)\int_{0}^t\|e(s)\|^2\mathrm{d}s$. The limit $t\!\to\!\infty$ of the last~in-\linebreak equality yields $\|e\|_{L^2}\leq\big({\tfrac{ V(\cdot)|_{t=0}}{\lambda_{\min}(K)}}\big)^{1/2}\!<\infty$.~That~is,~in~ad-\linebreak dition to $\dot{e}\in L^\infty$ per~\ref{Boundedness}),~\eqref{ErrorDyn}, one has $e\in L^2$.~Due to \cite{tao_simple_1997},\linebreak {this guarantees asymptotic tracking, i.e.\ $\underset{t\to\infty}{\lim}e(t)=0$}.
\end{enumerate} 
\vspace{-.75mm}
\subsection{Different data-acquisition methods to facilitate CL}\vspace{-.75mm}
The notable lack of estimation-error-convergence from the closed-loop properties~\ref{Stability})-\ref{AsymptoticTracking}) of the classical approach~with-\linebreak out PE motivates CL as initially proposed by the authors~of \cite{chowdhary_concurrent_2010,chowdhary_exponential_2014}: The dynamical behavior~\eqref{Plant} of the plant can be sampled and reorganized into a linear regression equation~(LRE)
\begin{align}
\label{basicLRE}
   Y(x(\tau))^\top\left(u(\tau)-\dot{x}(\tau)\right)&=-Y(x(\tau))^\top Y(x(\tau))\theta.
\end{align}
Therein, $\tau\in\mathcal{T}(t)\subseteq[\Delta t,t]$ is one of the $|\mathcal{T}(t)|<\infty$ historical data-acquisition instants, greater or equal some buffer constant $\Delta t\geq 0$, that are stored in a user-managed set $\mathcal{T}(\cdot)$.\linebreak The terms on the LHS of~\eqref{basicLRE} are, in principle, all available---provided, the state derivative~$\dot{x}$ is accessible---for being utilized as feedback in the update law to extract the information on the RHS of~\eqref{basicLRE} about the unknown parameters. Under availability of the state derivative, \cite{chowdhary_concurrent_2010} proposes to augment the update law~\eqref{UpdateLaw} with the sum of so-far-collected data~as~per
\begin{align}
\label{UpdateLawCL}
    \dot{\hat{\theta}}&=\Gamma\Big(Y(x)^\top e\ldots\\&\qquad\quad -k\textstyle\sum_{\tau \in\mathcal{T}(t)}Y(x(\tau))^\top\left(Y(x(\tau))\hat{\theta}+u(\tau)-\dot{x}(\tau)\right)\Big)\nonumber
\end{align}
with some gain $k(t)\geq 0$ that enables activating ($k(\cdot)>0$, constant) or deactivating ($k(\cdot)=0$) the CL-augmentation in a user-managed way. In view of the LRE~\eqref{basicLRE}, one then achieves
\begin{align}
\label{dotV_CL}
    \dot{V}(e,\tilde{\theta})&=-e^\top Ke-\tilde{\theta}^\top\Big(k\textstyle\sum_{\tau \in\mathcal{T}(t)}Y(x(\tau))^\top Y(x(\tau))\Big)\tilde{\theta}
\end{align}
from~\eqref{UpdateLawCL}. As the time derivative~\eqref{dotV_CL} admits $\dot{V}(\cdot)\leq-e^\top K e$, the closed loop of plant~\eqref{Plant} under adaptive control~\eqref{AdaptiveControl} with CL-augmented update law~\eqref{UpdateLawCL} does retain all of the properties of the classical approach listed in~\ref{Stability})-\ref{AsymptoticTracking}). Moreover, in addition to its negative definiteness w.r.t.\ the tracking error, the time derivative~\eqref{dotV_CL} of the Lyapunov function becomes negative definite also w.r.t.\ the estimation error---provided that a time instant $t^\star_\mathrm{CL}\in [\Delta t,\infty)$ exists during closed-loop operation, at which the following requirements are satisfied: At $t^\star_\mathrm{CL}$, the collected data as utilized by~\eqref{UpdateLawCL} is FE, i.e.\ $\sum_{\tau \in\mathcal{T}(t^\star_\mathrm{CL})}Y(x(\tau))^\top Y(x(\tau))\succ 0$, and the data collection is stopped thereafter while the CL-augmentation of~\eqref{UpdateLawCL} gets activated, i.e.\ $\mathcal{T}(t\geq t^\star_\mathrm{CL})=\mathcal{T}(t^\star_\mathrm{CL})$ and $k(t\geq t^\star_\mathrm{CL})=k(t^\star_\mathrm{CL})>0$. In a scenario where such an instant $t^\star_\mathrm{CL}$ exists, the full state $z$ converges toward zero exponentially afterwards; a fact that follows immediately from the consequence $\dot{V}(\cdot)|_{t\geq t^\star_\mathrm{CL}}\leq-\lambda^\star_\mathrm{CL} V(\cdot)|_{t\geq t^\star_\mathrm{CL}}$ of~\eqref{LyapFcn},~\eqref{dotV_CL} with exponential convergence rate
$$\lambda^\star_\mathrm{CL}=\tfrac{\min\{\lambda_{\min}(K), k(t^\star_\mathrm{CL})\lambda_{\min}(\sum_{\tau \in\mathcal{T}(t^\star_\mathrm{CL})}Y(x(\tau))\!^\top\! Y(x(\tau)))\}}{2^{-1}\scriptstyle\max\{1, \lambda_{\max}(\Gamma^{-1})\}}>0,$$
leading to exponential decay of $\|z\|^2\leq \tfrac{ 2}{\min\{1,\lambda_{\min}(\Gamma^{-1})\}}V(\cdot)$.

To enable CL~\eqref{UpdateLawCL} based on the LRE~\eqref{basicLRE} also for problems where the state derivative is not accessible,~\cite{chowdhary_exponential_2014} resorts to approximation via smoothing, whereas~\cite{kamalapurkar_concurrent_2017} proposes~estimation.

In~\cite{parikh_integral_2019}, a more advanced integral variation of CL, ICL,~is discovered that removes any reliance on the state derivative~$\dot{x}$: They construct another LRE, different from~\eqref{basicLRE} that allows extracting information about the unknown parameters from the closed-loop movement while circumventing $\dot{x}$ via appli-cation of the fundamental theorem of calculus (FTC).~Define
\begin{align}\label{FTCTrickery}
    \mathcal{X}(\tau)&=x(\tau)-x(\tau-\Delta t),\\
    \label{Rattenschwanz}
    \mathcal{U}(\tau)&=\textstyle\int_{\tau-\Delta t}^\tau u(s)\mathrm{d}s,\quad \mathcal{Y}(\tau)=\int_{\tau-\Delta t}^\tau Y(x(s))\mathrm{d}s.
\end{align}
Crucially, $\mathcal{X}$, due to FTC and absolute continuity of $x$,~admits $X(\tau)=\int_{\tau-\Delta t}^\tau \dot{x}(s)\mathrm{d}s$. Hence, by taking the integral over $[\tau-\Delta t,\tau]$ of the dynamical behavior in~\eqref{Plant}, one obtains the~LRE
\begin{align}
    \label{ICL_LRE}
    \mathcal{Y}(\tau)^\top\left(\mathcal{U}(\tau)-\mathcal{X}(\tau)\right)&=-\mathcal{Y}(\tau)^\top\mathcal{Y}(\tau)\theta,
\end{align}
whose LHS is computable via~\eqref{FTCTrickery},~\eqref{Rattenschwanz} without $\dot{x}$. Analogous to \eqref{UpdateLawCL}, the ICL-augmentation of~\eqref{UpdateLaw} with $\dot{x}$-free data~is
\begin{align}
    \label{UpdateLawICL}
    \dot{\hat{\theta}}&=\Gamma\Big(Y(x)^\top e-k\textstyle\sum_{\tau \in\mathcal{T}(t)}\mathcal{Y}(\tau)^\top\left(\mathcal{Y}(\tau)\hat{\theta}+\mathcal{U}(\tau)-\mathcal{X}(\tau)\right)\Big).
\end{align}
Then, adaptive control~\eqref{AdaptiveControl} with the ICL-update-law~\eqref{UpdateLawICL}, in view of \eqref{ICL_LRE}, prompts the Lyapunov function to evolve~with
\begin{align}
\label{dotV_ICL}
    \dot{V}(e,\tilde{\theta})&=-e^\top Ke-\tilde{\theta}^\top\Big(k\textstyle\sum_{\tau \in\mathcal{T}(t)}\textstyle\mathcal{Y}(\tau)^\top \mathcal{Y}(\tau)\Big)\tilde{\theta}.
\end{align}
Correspondingly, one must monitor the minimal~eigenvalue
\begin{align}
\label{lminICL}
\lambda^{\!\mathcal{YY}}(t)&=\lambda_{\min}\Big(\textstyle\sum_{\tau\in\mathcal{T}(t)}\mathcal{Y}(\tau)^\top\mathcal{Y}(\tau)\Big)
\end{align}in order to stay alert of any occurrence of FE in closed loop of the data as is utilized by the ICL-augmented~\eqref{UpdateLawICL}; which is evident should the eigenvalue in~\eqref{lminICL} become positive~at\linebreak some time instant $t^\star_\mathrm{ICL}\!\!\in\![\Delta t,\infty)$, i.e.\ $\lambda^{\!\mathcal{YY}}(t^\star_\mathrm{ICL})\!>\!0$.~In~a~scena-rio where such a $t^\star_\mathrm{ICL}$ exists, stopping the accumulation of data and activating ICL in~\eqref{UpdateLawICL}, i.e.\ $\mathcal{T}(t\geq t^\star_\mathrm{ICL})=\mathcal{T}(t^\star_\mathrm{ICL})$ and $k(t\geq t^\star_\mathrm{ICL})=k(t^\star_\mathrm{ICL})>0$,~guarantees (analogous to CL) the exponential convergence of the full state $z$ toward~zero.

In order to explain developments of CL beyond~ICL~as in \cite{pan_model_2016}, \cite{shen_regulation-triggered_2026}---based on yet another LRE, other than~\eqref{ICL_LRE}, that also facilitates $\dot{x}$-free information-extraction about the unknown parameters---the CL-survey \cite{wen_relaxation_2026} offers the following critique of ICL: The efficiency of the ICL-data-acquisition method \eqref{FTCTrickery}, \eqref{Rattenschwanz} may not be optimal since it first integrates the regressor $Y(x)$ to $\mathcal{Y}(\tau)$ before calculating the dyadic product $\mathcal{Y}(\tau)^\top \mathcal{Y}(\tau)$ that builds toward eventual positive definiteness of the matrix that determines FE. Specifically,~\cite{wen_relaxation_2026} points out that integrating first may cause the contributions $Y(x(s_i))$ of the regressor at distinct time instances $s_1\neq s_2$, $s_1,s_2\in[\tau-\Delta t,\tau]$ to cancel each other in~\eqref{Rattenschwanz}, thereby never entering the dyadic product that facilitates FE under ICL. The order of integration (to get rid of $\dot{x}$) and construction of dyadic products (to generate a symmetric, possibly positive definite, structure in $-\dot{V}(\cdot)$) can be reversed relative to~ICL~via
\begin{align}\label{JaAnscheinendIstDasSo}
    \mathcal{Z}(\tau)&=\textstyle\int_{x|_{[\tau-\Delta t,\tau]}}Y(\xi)^\top\mathrm{d}\xi,\\
    \label{ForDiscussionsSake}
    \mathcal{V}(\tau)&=\textstyle\int_{\tau-\Delta t}^\tau Y(x(s))^\top u(s)\mathrm{d}s,\\\label{HachJa}
    \mathcal{W}(\tau)&=\textstyle\int_{\tau-\Delta t}^\tau Y(x(s))^\top Y(x(s))\mathrm{d}s,
\end{align}
cf.\ \cite{wen_relaxation_2026}, as the matrix $\mathcal{W}$ therein integrates dyadic~products. The terms in~\eqref{JaAnscheinendIstDasSo}-\eqref{HachJa}  then satisfy the LRE with~$\dot{x}$-free~LHS
\begin{align}
\label{fullExcitation_LRE}
    \mathcal{V}(\tau)-\mathcal{Z}(\tau)&=-\mathcal{W}(\tau)\theta;
\end{align}
a fact that can be seen by integrating the reorganized version $Y(x)^\top(u-\dot{x})=-Y(x)^\top Y(x)\theta$ of~\eqref{Plant} over $[\tau-\Delta t,\tau]$ in view of $\mathcal{Z}(\tau)=\int_{\tau-\Delta t}^\tau Y(x(s))^\top \dot{x}(s)\mathrm{d}s$, which holds true due to\linebreak continuity of $Y\!$ and absolute continuity of $x$. However, forcing\linebreak the dyadic products into the integrand of $\mathcal{W}$ generally comes at the following cost: Extracting information about the~un-known $\theta$ from the RHS of~\eqref{fullExcitation_LRE} by computing its LHS,~while possible without $\dot{x}$, generally requires a recording of the oriented path of $x$ during $[\tau-\Delta t,\tau]$ as well as the calculation of the spatial tangential alongside this path in order to obtain the path integral $\mathcal{Z}(\tau)$---as opposed to the corresponding term of ICL, $\mathcal{X}(\tau)$, that can do only with the path's end points\linebreak $x(\tau-\Delta t)$, $x(\tau)$. Thus, the data-acquisition method~\eqref{JaAnscheinendIstDasSo}-\eqref{HachJa}, in general, is more demanding than data acquisition via~\eqref{FTCTrickery},~\eqref{Rattenschwanz} under ICL---except for cases in which the path integral in~\eqref{JaAnscheinendIstDasSo} may be simplified: E.g., if each column of the regressor happens to be a gradient field (see Assumption~\ref{GradientAssumption}) or by exploitation of a single-input structure as in~\cite{pan_model_2016}.
\begin{assumption}\label{GradientAssumption}
    The columns $Y_l:\mathbb{R}^n\to\mathbb{R}^n$ of $Y(\cdot)=(\begin{smallmatrix}Y_1(\cdot)&\ldots&Y_p(\cdot)\end{smallmatrix})$ are gradient fields, i.e.\ there exist potentials~$\Phi_l:\mathbb{R}^n\to\mathbb{R}$ with $Y_l(\xi)^\top\!=\tfrac{\partial}{\partial\xi}\Phi_l(\xi)$ for all $\xi\in\mathbb{R}^n$, $l\in\{1,\ldots,p\}$.
\end{assumption}
In case the regressor fulfills Assumption~\ref{GradientAssumption}, the path~integral $\mathcal{Z}(\tau)$ in~\eqref{JaAnscheinendIstDasSo} can be computed more easily as~per
\begin{align}
\label{PathInGradientFields}
    \mathcal{Z}(\tau)&=(\begin{smallmatrix}
        \Phi_1(x(\tau))-\Phi_1(x(\tau-\Delta t))&\ldots&\Phi_p(x(\tau))-\Phi_p(x(\tau-\Delta t))
    \end{smallmatrix})^\top
\end{align}
by using only the end points of $x|_{[\tau-\Delta t,\tau]}$. That is, Assumption~\ref{GradientAssumption} renders data acquisition to facilitate the update~law 
\begin{align}
    \label{UpdateLawNameWärNett}
    \dot{\hat{\theta}}&=\Gamma\Big(Y(x)^\top e-k\textstyle\sum_{\tau \in\mathcal{T}(t)}\textstyle\mathcal{W}(\tau)\hat{\theta}+\mathcal{V}(\tau)-\mathcal{Z}(\tau)\Big)
\end{align}
into being equally as demanding as it is for~\eqref{UpdateLawICL} under ICL.
Notably, adaptive control~\eqref{AdaptiveControl} with~\eqref{UpdateLawNameWärNett} serving as $\dot{x}$-free-data-\newpage\noindent augmented update law, via LRE~\eqref{fullExcitation_LRE}, does result in
\begin{align}
\label{dotV_NameWärNett}
    \dot{V}(e,\tilde{\theta})&=-e^\top Ke-\tilde{\theta}^\top\Big(k\textstyle\sum_{\tau \in\mathcal{T}(t)}\mathcal{W}(\tau)\Big)\tilde{\theta}.
\end{align}
Accordingly, the data in~\eqref{UpdateLawNameWärNett} is FE if the minimal eigenvalue
\begin{align}
\label{lminNameWärNett}
    \lambda^{\!\mathcal{W}}(t)&=\lambda_{\min}\Big(\textstyle\sum_{\tau\in\mathcal{T}(t)}\mathcal{W}(\tau)\Big)
\end{align}
happens to be positive. Provided, FE of the data occurs under adaptive control~\eqref{AdaptiveControl} with~\eqref{UpdateLawNameWärNett} of the plant~\eqref{Plant} at a time instant $t^\star\in[\Delta t,\infty)$, i.e.\ $\lambda^{\!\mathcal{W}}(t^\star)>0$, then imposing~$\mathcal{T}(t\!\geq\! t^\star)=\mathcal{T}(t^\star)$ and $k(t\!\geq\! t^\star)=k(t^\star)>0$, ensures (analogous to CL as well as to ICL) the exponential convergence of $z$~toward~zero.

\section{Results}\label{section:RESULTS}
Adaptive control~\eqref{AdaptiveControl} of the plant~\eqref{Plant}, as presented so far in Section \ref{section:PRELIMINARIES}, does not entail a-priori-verifiable conditions for closed-loop FE-occurrence when using any of the different data-augmented CL update laws~\eqref{UpdateLawCL}, \eqref{UpdateLawICL}, \eqref{UpdateLawNameWärNett}. Also, we seek to be independent of state-derivative accessibility. To address these points, we subsequently develop the reference-replay mechanism for the update laws~\eqref{UpdateLawICL}~and~\eqref{UpdateLawNameWärNett}, whose respective data can be acquired via~\eqref{FTCTrickery},~\eqref{Rattenschwanz}~and \eqref{JaAnscheinendIstDasSo}\footnote{If Assumption~\ref{GradientAssumption} holds, one should use the computationally lighter~\eqref{PathInGradientFields}.}-\eqref{HachJa} without the state derivative $\dot{x}$: Specifically, we propose to replay some pre-defined training segment $x_\mathrm{t}:[0,T_\mathrm{p}]\to\linebreak \mathbb{R}^n$ with $x_\mathrm{t}\in C^1([0,T_\mathrm{p}])$ of duration $T_\mathrm{p}\geq\Delta t$  
in the reference
\begin{align}
\label{Replay}
x_\mathrm{r}(t)&=\left\{\begin{matrix}x_\mathrm{t}\left(t-\left\lfloor t/T_\mathrm{p}\right\rfloor T_\mathrm{p}\right)&:\, t\in[0,T^\star)\\
x_ {\mathrm{d},T^\star}(t)&:\, t\in[T^\star,\infty)\end{matrix}\right.
\end{align}
until the first time instance of closed-loop FE-occurrence,
\begin{align}
\label{FirstFEOccurrence}
    \!\!T^\star&\!=\!\left\{\begin{matrix}
        \min\left\{t\geq \Delta t:\lambda^{\!\mathcal{YY}}(t)>0\right\}:&\!\!\!\!\text{$\dot{\hat{\theta}}$ as per \eqref{UpdateLawICL}}\\
        \min\left\{t\geq \Delta t:\lambda^{\!\mathcal{W}}(t)>0\right\}:&\!\!\!\!\text{$\dot{\hat{\theta}}$ as per~\eqref{UpdateLawNameWärNett}}
    \end{matrix}\right.,
\end{align}
after which reference replay as per~\eqref{Replay} permits to track any desired reference trajectory
$x_{\mathrm{d},T^\star}:[T^\star,\infty)\to\mathbb{R}^n$ that fulfills the trackability conditions $x_{\mathrm{d},T^\star}\in C^1([T^\star,\infty))\cap L^\infty([T^\star,\infty))$ and $\dot{x}_{\mathrm{d},T^\star}\in L^\infty([T^\star,\infty))$. In order to guarantee the trackability of the reference also before FE occurs in closed loop, i.e. during the replay phase $[0,T^\star)$ of~\eqref{Replay}, we require the training segment $x_\mathrm{t}$ to fulfill the boundary conditions (BCs) $x_\mathrm{t}(T_\mathrm{p})=x_\mathrm{t}(0)$ and $\lim_{t\nearrow T_\mathrm{p}}\dot{x}_\mathrm{t}(t)=\lim_{t\searrow 0}\dot{x}_\mathrm{t}(t)$.  Together with these BCs, $x_\mathrm{t}\in C^1([0,T_\mathrm{p}])$ then ensures $x_\mathrm{r}([0,T^\star))\in C^1([0,T^\star))$~in addition to the consequence $x_\mathrm{r},\dot{x}_\mathrm{r}\!\in\! L^\infty([0,T^\star))$ resulting~from \eqref{Replay} and the boundedness inherent to $x_\mathrm{t}\in C^1([0,T_\mathrm{p}])$ and~$\dot{x}_\mathrm{t}$. Thereby, we have trackability of $x_\mathrm{r}$ before and after $T^\star$ under usage of the reference-replay mechanism~\eqref{Replay}. However, note that~\eqref{Replay} allows non-trackability of $x_\mathrm{r}$ at $t=T^\star$, i.e.\ we permit arbitrary initial conditions $x_{\mathrm{d},T^\star}(T^\star),\dot{x}_{\mathrm{d},T^\star}(T^\star)\in\mathbb{R}^n$ of the desired reference that are not in smooth agreement with $x_\mathrm{r}$ for $t\nearrow T^\star$; which we treat as a reinitialization of the respective closed loop (for either update law,~\eqref{UpdateLawICL}~or \eqref{UpdateLawNameWärNett}) at $T{^\star}$.

\subsection{A-priori-verifiable conditions for the occurrence of FE}
Here, we present a-priori-verifiable conditions to guarantee FE-occurrence at $T^\star<\infty$ under reference replay~\eqref{Replay} when using adaptation~\eqref{AdaptiveControl} with either~\eqref{UpdateLawICL} or~\eqref{UpdateLawNameWärNett} as update~law.

Let $\mathcal{T}_\mathrm{t}\subset[\Delta t\!, T_\mathrm{p}]$ be a pre-defined set of $|\mathcal{T}_\mathrm{t}|<\infty$ time~instan-ces at which data is acquired from the training segment~via
\begin{align}
    \mathcal{Y}_\mathrm{t}(\tau_\mathrm{t})&=\textstyle\int_{\tau_\mathrm{t}-\Delta t}^{\tau_\mathrm{t}} Y(x_\mathrm{t}(s))\mathrm{d}s,\\
    \mathcal{W}_\mathrm{t}(\tau_\mathrm{t})&=\textstyle\int_{\tau_\mathrm{t}-\Delta t}^{\tau_\mathrm{t}} Y(x_\mathrm{t}(s))^\top Y(x_\mathrm{t}(s))\mathrm{d}s \text{ for }\tau_\mathrm{t}\in\mathcal{T}_\mathrm{t}.
\end{align}
As both the training segment $x_\mathrm{t}$ and the set $\mathcal{T}_\mathrm{t}$ are available in advance, one then can compute the minimal~eigenvalues
\begin{align}
	\label{TrainingICL}\lambda_{\,\mathrm{t}}^{\!\mathcal{YY}}&=\lambda_{\min}\Big(\textstyle\sum_{\tau_\mathrm{t}\in\mathcal{T}_\mathrm{t}}\mathcal{Y}_\mathrm{t}(\tau_\mathrm{t})^\top \mathcal{Y}_\mathrm{t}(\tau_\mathrm{t})\Big),\\
    \label{TrainingNameWärNett}\lambda_{\,\mathrm{t}}^{\!\mathcal{W}}&=\lambda_{\min}\Big(\textstyle\sum_{\tau_\mathrm{t}\in\mathcal{T}_\mathrm{t}}\mathcal{W}_\mathrm{t}(\tau_\mathrm{t})\Big)
\end{align}
offline. This facilitates a-priori-verifiability of the excitation requirement for the training segment; which is met if $\lambda_{\,\mathrm{t}}^{\!\mathcal{YY}}$ or $\lambda_{\,\mathrm{t}}^{\!\mathcal{W}}$ has a positive sign---depending on the user's choice of update law: $\lambda_{\,\mathrm{t}}^{\!\mathcal{YY}}$ for~\eqref{UpdateLawICL} and $\lambda_{\,\mathrm{t}}^{\!\mathcal{W}}$ for~\eqref{UpdateLawNameWärNett}. To ensure~that the full training-excitation (encoded by $\lambda_{\,\mathrm{t}}^{\!\mathcal{YY}}$, $\lambda_{\,\mathrm{t}}^{\!\mathcal{W}}$) is injected into~the data collected on the closed-loop trajectory $x=x_\mathrm{r}+e$, hence guaranteeing FE once $e$  becomes sufficiently small, we moreover require the data-acquisition instants $\mathcal{T}(\cdot)$ to satisfy
\begin{align}
\label{DataReplay}
\!\!\left\{\!\!\begin{matrix}\mathcal{T}(t)\supseteq\left(\mathcal{T}_\mathrm{t}\oplus\left\{\left\llbracket\tfrac{t}{T_\mathrm{p}}\right\rrbracket T_\mathrm{p}\right\}\right)\!\!\!\!\!\!\!\!&:\, t\in\left\{T_{\max}\!+\!lT_\mathrm{p}\!:l\!\in\!\mathbb{N}_0\right\}\cap[0,\!T^\star)\\
\raggedleft\mathcal{T}(t)=\mathcal{T}(T^\star)\!\!\!\!\!\!\!\!&:\, t\in[T^\star,\infty)\end{matrix}\right.
\end{align}
Therein, $T_{\max}\!=\max\!\mathcal{T}_\mathrm{t}$ is the time instant at which~data~acqui-sition from the training segment is concluded. Crucially, the first line of \eqref{DataReplay} requires including suitably translated training\linebreak data-acquisition instants in $\mathcal{T}(\cdot)$ for each replay cycle at the specific instant when the full excitation of the training can be referenced,~e.g.\ at $T_{\max}+lT_\mathrm{p}$ for the $l$th cycle. We note~that
\begin{align}
\label{ExplicitDataReplay}
\!\!\!\!\mathcal{T}(t)&=\!\left\{\!\begin{matrix}\left(\!\mathcal{T}_\mathrm{t}\cap\left[0,t\!-\!\left\llbracket\!\tfrac{t}{T_\mathrm{p}}\!\right\rrbracket T_\mathrm{p}\right]\!\right)\oplus \left\{\left\llbracket\!\tfrac{t}{T_\mathrm{p}}\!\right\rrbracket T_\mathrm{p}\right\}\!\!\!\!\!\!\!\!\!&:\, t\in[0,T^\star)\\
\mathcal{T}(T^\star)&:\, t\in[T^\star,\infty)\end{matrix}\right.\!\!\!
\end{align}
is an implementable example for a $\!\mathcal{T}\!(\cdot)\!$ that fulfills~the~replay-\linebreak requirement~\eqref{DataReplay} while also respecting causality and the storage limitations from Section~\ref{section:PRELIMINARIES}, i.e.\ $\mathcal{T}(t)\subseteq[\Delta t, t]$ and $|\mathcal{T}(t)|\leq|\mathcal{T}_\mathrm{t}|<\infty$ for all~$\!t\geq 0$; all of which is readily~apparent.

\begin{remark}
    We use our custom-defined function $\llbracket\cdot\rrbracket$ instead of the more common floor-function $\lfloor\cdot\rfloor$ in \eqref{DataReplay}, \eqref{ExplicitDataReplay} for the sole reason to allow treating the edge case $T_{\max}=T_\mathrm{p}$ with the subsequent proofs of Theorem~\ref{Theorem:ICL},~\ref{Theorem:NameWärNett}  and~\ref{theorem:EmployingExpertKnowledge}, as they rely on the validity of $\llbracket (T_{\max}+lT_\mathrm{p})/T_\mathrm{p}\rrbracket=l$ for all $l\in\mathbb{N}_0=\mathbb{N}\cup\{0\}$.\linebreak Thus, one may replace $\llbracket\cdot\rrbracket$ in~\eqref{DataReplay}, $\!$\eqref{ExplicitDataReplay} by $\lfloor\cdot\rfloor$  if $T_{\max}\!\in\![\Delta t\!,T_\mathrm{p})$.
\end{remark}

First, we show our central result for ICL under reference replay (i.e. with~\eqref{UpdateLawICL} as $\dot{x}$-free-data-augmented update law):
\begin{theorem}\label{Theorem:ICL}
Adaptive control~\eqref{AdaptiveControl} of~\eqref{Plant} under ICL~\eqref{UpdateLawICL} and reference replay~\eqref{Replay} guarantees the boundedness of the full state $z=(e^\top,\tilde{\theta}^\top)^\top$ of the closed loop as well as the asymptotic convergence of $e$ toward zero. Moreover, if the training segment passes the offline excitation test $\lambda_{\,\mathrm{t}}^{\!\mathcal{YY}}>0$ and the data-acquisition instants fulfill the replay-requirement~\eqref{DataReplay}, then FE does occur at $T^\star<\infty$ such that activating the ICL-augmentation in~\eqref{UpdateLawICL} via $k(t\geq T^\star)=k(T^\star)>0$ subsequently leads to exponential convergence of the full state with $\forall\,t\geq\allowbreak T^\star:\|z(t)\|\leq V(z(T^\star))^{1/2}M^\star e^{c^\star_\mathrm{ICL}(T^\star-t)}$, where $M^\star,c^\star_\mathrm{ICL}>0$.\linebreak Notably, the convergence rate $c^\star_\mathrm{ICL}$ (as is evident by its struc-\linebreak ture in the proof) is arbitrarily large by choice of~$K$,~$k(T^\star)$.
\end{theorem}

\begin{proof}
The claims are validated in several stages:

1a) \textit{Boundedness.} Since the time derivative~\eqref{dotV_ICL} of $V(\cdot)$ is nonpositive at all time instants except $T^\star$ where $\dot{V}$ may not exist, the argument of \ref{Boundedness}) still delivers boundedness of~$z$.

1b) \textit{Asymptotic tracking.} First, we consider the scenario $T^\star=\infty$. Then, $x_\mathrm{r}$ as per~\eqref{Replay} is an element of $C^1\cap L^\infty$ due to the BCs imposed on $x_\mathrm{t}$. This reference furthermore satisfies $\dot{x}_\mathrm{r}\in L^\infty$ as it is a replayed version of $x_\mathrm{t}\in C^1([0,T_\mathrm{p}])$ which, inherently, has a bounded time derivative. Hence, $x_\mathrm{r}$ is trackable and thus, we obtain $\lim_{t\to\infty}e(t)=0$  by arguing as in \ref{AsymptoticTracking}). Second, we consider the only remaining scenario $T^\star<\infty$. The argumentative opening of \ref{AsymptoticTracking})---when replacing the initial time $0$ by $T^\star$---yields $\int_{T^\star}^\infty \|e(s)\|^2\mathrm{d}s\leq \tfrac{V(\cdot)|_{t=T^\star}}{\lambda_{\min}(K)}<\infty$, i.e.\ $e\in L^2([T^\star,\infty))$. Furthermore, the continuity of $Y$ gives $Y(x)\in L^\infty([T^\star,\infty))$ since $x|_{t\geq T^\star}=e+x_{\mathrm{d},T^\star}\in L^\infty([T^\star,\infty))$ due to $e$ being bounded by virtue of 1a). Thus, the tracking-error dynamics~\eqref{ErrorDyn} ensure $\dot{e}\in L^\infty([T^\star,\infty))$ as 1a) also guarantees the boundedness of $\tilde{\theta}$. Via~\cite{tao_simple_1997}, this leads to $\lim_{t\to\infty}e(t)=0$. 

1c) \textit{FE-occurrence.} For the remainder of this proof, let \eqref{DataReplay}$\!$ apply and $\lambda_{\,\mathrm{t}}^{\!\mathcal{YY}}\!>\!0$. We aim to validate $T^\star\!<\!\infty$~by contradict-\linebreak ion and thus, assume $T^\star\!=\!\infty$. Define $\mathcal{S}_{\Delta t}\!:\![\Delta t\!,\infty)\!\to\![0,\infty)$~with
\begin{align}
\label{SupremeMaxSing}
    {\scriptstyle\mathcal{S}_{\Delta t}(t)}&={\scriptstyle\underset{s\in[t-\Delta t,t]}{\sup}\left\{\sigma_{\max}\left({Y}(x_\mathrm{r}(s))-Y(x(s))\right)\right\}},
\end{align}
as it appears in the following: Whenever $t=T_{\max}+lT_\mathrm{p}$ for a $l\in\mathbb{N}_0$, then~\eqref{Replay},~\eqref{DataReplay} in view of $t<T^\star\!$  yield the lower~bound
\begin{align}
    &{\scriptstyle\vartheta^\top\big(\underset{\tau\in\mathcal{T}(t)}{\sum}\mathcal{Y}(\tau)^\top\mathcal{Y}(\tau)\big)\vartheta}\overset{\scriptscriptstyle\eqref{DataReplay}}{\geq}{\scriptstyle\underset{\tau\in(\mathcal{T}_\mathrm{t}\oplus\{\llbracket t/T_\mathrm{p}\rrbracket T_\mathrm{p}\})}{\sum}\vartheta^\top\mathcal{Y}(\tau)^\top\mathcal{Y}(\tau)\vartheta}\nonumber\allowdisplaybreaks\\
    &={\scriptstyle\underset{\tau_\mathrm{t}\in \mathcal{T}_\mathrm{t}}{\sum}\vartheta^\top\mathcal{Y}(\tau_\mathrm{t}+lT_\mathrm{p})^\top\mathcal{Y}(\tau_\mathrm{t}+lT_\mathrm{p})\vartheta}\nonumber\allowdisplaybreaks\\
    &={\scriptstyle\underset{\tau_\mathrm{t}\in \mathcal{T}_\mathrm{t}}{\sum}\vartheta^\top\left(\mathcal{Y}_\mathrm{t}(\tau_\mathrm{t})+\mathcal{Y}(\tau_\mathrm{t}+lT_\mathrm{p})-\mathcal{Y}_\mathrm{t}(\tau_\mathrm{t})\right)^\top\left(\mathcal{Y}_\mathrm{t}(\tau_\mathrm{t})+\mathcal{Y}(\tau_\mathrm{t}+lT_\mathrm{p})-\mathcal{Y}_\mathrm{t}(\tau_\mathrm{t})\right)\vartheta}\nonumber\allowdisplaybreaks\\
    &\geq {\scriptstyle\underset{\tau_\mathrm{t}\in \mathcal{T}_\mathrm{t}}{\sum}\vartheta^\top\mathcal{Y}_\mathrm{t}(\tau_\mathrm{t})^\top \mathcal{Y}_\mathrm{t}(\tau_\mathrm{t})\vartheta+ 2 \vartheta^\top \mathcal{Y}_\mathrm{t}(\tau_\mathrm{t})^\top\left(\mathcal{Y}(\tau_\mathrm{t}+lT_\mathrm{p})-\mathcal{Y}_\mathrm{t}(\tau_\mathrm{t})\right)\vartheta}\nonumber\allowdisplaybreaks\\
    &\geq {\scriptstyle\lambda_{\,\mathrm{t}}^{\!\mathcal{YY}} - 2 \underset{\tau_\mathrm{t}\in \mathcal{T}_\mathrm{t}}{\sum}\left\|\mathcal{Y}_\mathrm{t}(\tau_\mathrm{t})\vartheta\right\|\left\|\left(\mathcal{Y}_\mathrm{t}(\tau_\mathrm{t})-\mathcal{Y}(\tau_\mathrm{t}+lT_\mathrm{p})\right)\vartheta\right\|}\nonumber\allowdisplaybreaks\\
    &\geq {\scriptstyle\lambda_{\,\mathrm{t}}^{\!\mathcal{YY}} - 2 \underset{\tau_\mathrm{t}\in \mathcal{T}_\mathrm{t}}{\sum}\sigma_{\max}(\mathcal{Y}_\mathrm{t}(\tau_\mathrm{t}))\left\|\int_{\tau_\mathrm{t}-\Delta t}^{\tau_\mathrm{t}}Y(x_\mathrm{t}(s))\vartheta\mathrm{d}s-\int_{\tau_\mathrm{t}+lT_\mathrm{p}-\Delta t}^{\tau_\mathrm{t}+lT_\mathrm{p}}Y(x(s))\vartheta\mathrm{d}s\right\|}\nonumber\allowdisplaybreaks\\
    \label{AufhängerTheoremICL}
    &\overset{\scriptscriptstyle\eqref{Replay}}{=} {\scriptstyle\lambda_{\,\mathrm{t}}^{\!\mathcal{YY}} - 2 \underset{\tau_\mathrm{t}\in \mathcal{T}_\mathrm{t}}{\sum}\sigma_{\max}(\mathcal{Y}_\mathrm{t}(\tau_\mathrm{t}))\left\|\int_{\tau_\mathrm{t}+lT_\mathrm{p}-\Delta t}^{\tau_\mathrm{t}+lT_\mathrm{p}}\left(Y(x_\mathrm{r}(s))-Y(x(s))\right)\vartheta\mathrm{d}s\right\|}\allowdisplaybreaks\\
    &\geq {\scriptstyle\lambda_{\,\mathrm{t}}^{\!\mathcal{YY}} - 2\Delta t \underset{\tau_\mathrm{t}\in \mathcal{T}_\mathrm{t}}{\sum}\sigma_{\max}(\mathcal{Y}_\mathrm{t}(\tau_\mathrm{t}))\underset{ s\in[\tau_t+lT_\mathrm{p}-\Delta t,\tau_t+lT_\mathrm{p}]}{\!\!\!\!\!\!\!\!\!\!\!\!\!\!\!\!\sup}\!\!\!\!\!\!\!\!\!\!\!\!\!\!\!\!\!\!\!\!\!\!\!\!\!\!\left\{\sigma_{\max}\left({Y}(x_\mathrm{r}(s))-Y(x(s))\right)\right\}}\nonumber\allowdisplaybreaks\\
    &\geq{\scriptstyle\lambda_{\,\mathrm{t}}^{\!\mathcal{YY}} - 2\Delta t \big(\underset{\tau_\mathrm{t}\in \mathcal{T}_\mathrm{t}}{\sum}\sigma_{\max}(\mathcal{Y}_\mathrm{t}(\tau_\mathrm{t}))\big)\underset{\tau_\mathrm{t}\in\mathcal{T}_\mathrm{t}}{\max}\left\{\mathcal{S}_{\Delta t}(\tau_\mathrm{t}+lT_\mathrm{p})\right\}}
\end{align}
for any $\vartheta\in\mathbb{R}^p$ with $\|\vartheta\|=1$. Thus, we are sure~that
\begin{align}\label{LowerBoundMinEV_YY}
	\!\!\!{\scriptstyle\lambda^{\!\mathcal{YY}}(T_{\max}+lT_\mathrm{p})}\geq&{\scriptstyle\lambda_{\,\mathrm{t}}^{\!\mathcal{YY}} - 2\Delta t \big(\underset{\tau_\mathrm{t}\in \mathcal{T}_\mathrm{t}}{\sum}\sigma_{\max}(\mathcal{Y}_\mathrm{t}(\tau_\mathrm{t}))\big)\underset{\tau_\mathrm{t}\in\mathcal{T}_\mathrm{t}}{\max}\left\{\mathcal{S}_{\Delta t}(\tau_\mathrm{t}+lT_\mathrm{p})\right\}}\!\!\!
\end{align}
for all $l\in\mathbb{N}_0$. Consider $T_{\min}=\min\mathcal{T}_\mathrm{t}$ and let $\epsilon=(1-\mu)\lambda_{\,\mathrm{t}}^{\!\mathcal{YY}}/(2\Delta t\sum_{\tau_\mathrm{t}\in\mathcal{T}_\mathrm{t}}\sigma_{\max}(\mathcal{Y}_\mathrm{t}(\tau_\mathrm{t})))>0$ for a fixed $\mu\in(0,1)$. Then, in view of 1b) and Lemma~\ref{ReplayBisDerArztKommt} from the Appendix, the term on the RHS of~\eqref{LowerBoundMinEV_YY} is greater than or equal to $\mu \lambda_{\,\mathrm{t}}^{\!\mathcal{YY}}$ for all $l\geq (T(\epsilon)-T_{\min})/T_\mathrm{p}\geq 0$ (because the arguments of the different appearances of $\mathrm{S}_{\Delta t}$ on the RHS of~\eqref{LowerBoundMinEV_YY} then all fulfill $\tau_\mathrm{t}+lT_\mathrm{p}\geq T(\epsilon)$). Accordingly, we obtain the contradiction $\lambda^{\!\mathcal{YY}}(T_{\max}+lT_\mathrm{p})\geq\mu\lambda_{\,\mathrm{t}}^{\!\mathcal{YY}}>0$ to $T^\star=\infty$ at $l=\lceil(T(\epsilon)-T_{\min})/T_\mathrm{p}\rceil\in\mathbb{N}_0$; which in turn proves that FE will occur at some $T^\star<\infty$ during closed-loop operation.

1d) \textit{Exponential convergence subsequent to $T^\star$.} Provided $k(t\geq T^\star)=k(T^\star)$, we obtain $\dot{V}(\cdot)|_{t\geq T^\star}\leq-\lambda^\star_\mathrm{ICL} V(\cdot)|_{t\geq T^\star}$ with $\lambda^\star_\mathrm{ICL}=\tfrac{\min\{\lambda_{\min}(K), k(T^\star)\lambda^{\!\mathcal{YY}}(T^\star)\}}{2^{-1}\max\{1, \lambda_{\max}(\Gamma^{-1})\}}$ by means of~\eqref{LyapFcn},~\eqref{dotV_ICL} since $\mathcal{T}(t\geq T^\star)=\mathcal{T}(T^\star)$ holds by virtue of~\eqref{DataReplay}. Together with the consequence $\|z\|^2\leq \tfrac{2}{\min\{1,\lambda_{\min}(\Gamma^{-1})\}}V(\cdot)$ of~\eqref{LyapFcn}, this results in the claimed exponential convergence where the constants are  $M^\star={\scriptstyle\sqrt{2/\min\{1,\lambda_{\min}(\Gamma^{-1})\}}}$ and $c_\mathrm{ICL}^\star=\lambda^\star_\mathrm{ICL}/2$.
\end{proof}

Second, our central result also holds when using the $\dot{x}$-free-data-augmented update law \eqref{UpdateLawNameWärNett} under reference replay:

\begin{theorem}\label{Theorem:NameWärNett}
Let the plant~\eqref{Plant} be in closed loop under reference replay~\eqref{Replay} with adaptive control~\eqref{AdaptiveControl} where~\eqref{UpdateLawNameWärNett} serves as update law. Then, the full state $z$ remains bounded and $e$ converges asymptotically toward zero. Moreover, if the training segment passes the offline excitation test $\lambda_{\,\mathrm{t}}^{\!\mathcal{W}}>0$ and the data-acquisition instants fulfill the replay-requirement~\eqref{DataReplay}, then FE does occur at $T^\star<\infty$ such that activating the data-augmentation in~\eqref{UpdateLawNameWärNett} via $k(t\geq T^\star)=k(T^\star)>0$ subsequently leads to exponential convergence of the full state with $\forall\,t\geq T^\star\!:\|z(t)\|\!\leq\! V(z(T^\star))^{1/2}M^\star e^{c^\star(T^\star-t)}\!$; where $M^\star\!\!,c^\star\!>\!0$ and the\linebreak convergence rate $c^\star\!$ is arbitrarily large by choice of~$K$,~$k(T^\star)$.
\end{theorem}
\begin{proof}
2a) \textit{Boundedness.} $\dot{V}(\cdot)$ as per~\eqref{dotV_NameWärNett} being nonpositive for all times (except for $t=T^\star\!$~in case of $\lim_{t\nearrow T^\star}e(t)\neq e(T^\star)$) yields boundedness of $V(\cdot)$ at all times, including $T^\star$. Thus, we obtain the claimed boundedness of $z$ from the (in \ref{Boundedness}) established) boundedness of the sublevel set $\Omega_{\sup_{t\geq 0}\{V(\cdot)|_{t}\}}$ of $V(\cdot)$; i.e.\ the argument, so far, is analogous to that~of~1a).~Furthermore, we define the terms $M_\mathrm{r}=\sup_{t\geq 0}\{\sigma_{\max}(Y(x_\mathrm{r}(t)))\}$,\linebreak $M=\sup_{t\geq 0}\{\sigma_{\max}(Y(x(t)))\}$, as they are needed below in~2c). We proceed by proving boundedness: Firstly, $M_\mathrm{r}<\infty$~follows from $x_\mathrm{r}$ as per \eqref{Replay} being bounded and $Y$ being continuous. Secondly, $x=e+x_\mathrm{r}$ is bounded since $z$ is bounded, such that $M<\infty$ then follows (again) from $Y$ being continuous. 

2b) \textit{Asymptotic tracking.} In either scenario, $T^\star=\infty$ or $T^\star<\infty$, we obtain the claimed $\lim_{t\to\infty}e(t)=0$ analogously to 1b) as the required boundedness of $e$ and $\tilde{\theta}$ holds by 2a).

2c) \textit{FE-occurrence.} Until the end of this proof, let \eqref{DataReplay}~ap-ply and $\lambda_{\,\mathrm{t}}^{\!\mathcal{W}}\!>0$. To show the claim $T^\star<\infty$ by contradiction, we assume $T^\star\!=\!\infty$. For any $t=T_{\max}+lT_\mathrm{p}$ with $l\in\mathbb{N}_0$, then
\begin{align}
    &{\scriptstyle\vartheta^\top\big(\underset{\tau\in\mathcal{T}(t)}{\sum}\!\!\!\mathcal{W}(\tau)\big)\vartheta}\overset{\scriptscriptstyle\eqref{DataReplay}}{\geq}{\scriptstyle\!\!\!\!\!\!\underset{\tau\in(\mathcal{T}_\mathrm{t}\oplus\{\llbracket t/T_\mathrm{p}\rrbracket T_\mathrm{p}\})}{\sum}\!\!\!\!\!\!\!\!\!\!\!\!\vartheta^\top\mathcal{W}(\tau)\vartheta}={\scriptstyle\underset{\tau_\mathrm{t}\in\mathcal{T}_\mathrm{t}}{\sum}\!\!\!\vartheta^\top\mathcal{W}(\tau_\mathrm{t}+lT_\mathrm{p})\vartheta}\allowdisplaybreaks\nonumber\\
    &\mkern-9mu= {\scriptstyle\underset{\tau_\mathrm{t}\in\mathcal{T}_\mathrm{t}}{\sum}\vartheta^\top\left(\mathcal{W}_\mathrm{t}(\tau_\mathrm{t})+\mathcal{W}(\tau_\mathrm{t}+lT_\mathrm{p})-\mathcal{W}_\mathrm{t}(\tau_\mathrm{t})\right)\vartheta}\allowdisplaybreaks\nonumber\\
    &\mkern-9mu\geq {\scriptstyle\lambda_{\,\mathrm{t}}^{\!\mathcal{W}}\!\!-\!\!\underset{\tau_\mathrm{t}\in\mathcal{T}_\mathrm{t}}{\sum}\left|\vartheta^{\!\top}\!\left(\int_{\tau_\mathrm{t}+lT_\mathrm{p}-\Delta t}^{\tau_\mathrm{t}+lT_\mathrm{p}}Y(x(s))^\top Y(x(s))\mathrm{d}s-\int_{\tau_\mathrm{t}-\Delta t}^{\tau_\mathrm{t}}Y(x_\mathrm{t}(s))^\top Y(x_\mathrm{t}(s))\mathrm{d}s\right)\!\vartheta \right|}\allowdisplaybreaks\nonumber\\
    &\mkern-9mu\overset{\scriptscriptstyle\eqref{Replay}}{\geq} {\scriptstyle\lambda_{\,\mathrm{t}}^{\!\mathcal{W}}\!-\!\underset{\tau_\mathrm{t}\in\mathcal{T}_\mathrm{t}}{\sum}\left|\vartheta^\top\left(\int_{\tau_\mathrm{t}+lT_\mathrm{p}-\Delta t}^{\tau_\mathrm{t}+lT_\mathrm{p}}Y(x(s))^\top Y(x(s))-Y(x_\mathrm{r}(s))^\top Y(x_\mathrm{r}(s))\mathrm{d}s\right)\vartheta \right|}\allowdisplaybreaks\nonumber\\
    &\mkern-9mu= {\scriptstyle\lambda_{\,\mathrm{t}}^{\!\mathcal{W}}\!-\!\underset{\tau_\mathrm{t}\in\mathcal{T}_\mathrm{t}}{\sum}\left|\int_{\tau_\mathrm{t}+lT_\mathrm{p}-\Delta t}^{\tau_\mathrm{t}+lT_\mathrm{p}}\vartheta^\top\left(Y(x_\mathrm{r}(s))+Y(x(s))\right)^\top\left(Y(x_\mathrm{r}(s))-Y(x(s))\right)\vartheta \mathrm{d}s \right|}\allowdisplaybreaks\nonumber\\
    \label{AufhängerTheoremNameWärNett}
    &\mkern-9mu\geq {\scriptstyle\lambda_{\,\mathrm{t}}^{\!\mathcal{W}}\!-\!\underset{\tau_\mathrm{t}\in\mathcal{T}_\mathrm{t}}{\sum}\int_{\tau_\mathrm{t}+lT_\mathrm{p}-\Delta t}^{\tau_\mathrm{t}+lT_\mathrm{p}}\left\|\left(Y(x_\mathrm{r}(s))+Y(x(s))\right)\vartheta\right\|\left\|\left(Y(x_\mathrm{r}(s))-Y(x(s))\right)\vartheta\right\| \mathrm{d}s }\!\!\!\allowdisplaybreaks\\
    &\mkern-9mu\geq {\scriptstyle\lambda_{\,\mathrm{t}}^{\!\mathcal{W}}-\Delta t(M_\mathrm{r}+M) \underset{\tau_\mathrm{t}\in \mathcal{T}_\mathrm{t}}{\sum}\quad\underset{ s\in[\tau_t+lT_\mathrm{p}-\Delta t,\tau_t+lT_\mathrm{p}]}{\!\!\!\!\!\!\!\!\!\!\!\!\!\!\!\!\sup}\!\!\!\!\!\!\!\!\!\!\!\!\!\!\!\!\!\!\!\!\!\!\!\!\!\!\left\{\sigma_{\max}\left({Y}(x_\mathrm{r}(s))-Y(x(s))\right)\right\}}\allowdisplaybreaks\nonumber\\
    &\mkern-9mu\geq {\scriptstyle\lambda_{\,\mathrm{t}}^{\!\mathcal{W}}-\Delta t (M_\mathrm{r}+M) |\mathcal{T}_\mathrm{t}|\underset{\tau_\mathrm{t}\in\mathcal{T}_\mathrm{t}}{\max}\left\{\mathcal{S}_{\Delta t}(\tau_\mathrm{t}+lT_\mathrm{p})\right\}}
\end{align}
holds for all $\vartheta\in\mathbb{R}^p$ with $\|\vartheta\|=1$ due to $t<T^\star$. This ensures
\begin{align}\label{LowerBoundMinEV_W}
	{\scriptstyle \lambda^{\!\mathcal{W}}(T_{\max}+lT_\mathrm{p}})&\geq{\scriptstyle\lambda_{\,\mathrm{t}}^{\!\mathcal{W}}-\Delta t (M_\mathrm{r}+M) |\mathcal{T}_\mathrm{t}|\underset{\tau_\mathrm{t}\in\mathcal{T}_\mathrm{t}}{\max}\left\{\mathcal{S}_{\Delta t}(\tau_\mathrm{t}+lT_\mathrm{p})\right\}}
\end{align}
for all $l\in\mathbb{N}_0$. Let $\epsilon=(1-\nu)\lambda_{\,\mathrm{t}}^{\!\mathcal{W}}/(\Delta t(M_\mathrm{r}+M)|\mathcal{T}_\mathrm{t}|)$ in which $\nu\in(0,1)$ is some fixed value. We have $\epsilon>0$ since $M_\mathrm{r},M<\infty$ as per 2a). Applying Lemma~\ref{ReplayBisDerArztKommt} from the Appendix in view of 2b) then yields existence of $T(\epsilon)\geq T_\mathrm{p}$ such that the term on the RHS of~\eqref{LowerBoundMinEV_W} is greater than or equal to $\nu\lambda_{\,\mathrm{t}}^{\!\mathcal{W}}$ for all $l\geq (T(\epsilon)-T_{\min})/T_\mathrm{p}\geq 0$ with $T_{\min}=\min\mathcal{T}_\mathrm{t}$. Analogous to the reasoning in 1c), we get $\lambda^{\!\mathcal{W}}(T_{\max}+\lceil (T(\epsilon)-T_{\min})/T_\mathrm{p}\rceil T_\mathrm{p})\geq \nu\lambda_{\,\mathrm{t}}^{\!\mathcal{W}}>0$ in contradiction to $T^\star=\infty$; thus establishing closed-loop FE-occurrence at some $T^\star<\infty$.

2d) \textit{Exponential convergence subsequent to $T^\star$.} Provided $k(t\geq T^\star)=k(T^\star)$, we obtain the exponential convergence of $z$ toward zero as claimed by arguing as in 1d)---but for the time derivative $\dot{V}(\cdot)$ as per~\eqref{dotV_NameWärNett} (instead of~\eqref{dotV_ICL})---with the constants $M^\star$ as in 1d) and $c^\star=\tfrac{\min\{\lambda_{\min}(K), k(T^\star)\lambda^{\!\mathcal{W}}(T^\star)\}}{\max\{1, \lambda_{\max}(\Gamma^{-1})\}}$.
\end{proof}

\subsection{Upper bound on $T^\star$ under local Lipschitz continuity of the regressor (i.e., requiring additional smoothness)}
Here, we derive (conservative) upper bounds on the first time instants  $T^\star<\infty$ at which the a-priori-verifiable conditions in Theorems~\ref{Theorem:ICL},~\ref{Theorem:NameWärNett} guarantee FE in closed loop under \eqref{Replay}, \eqref{DataReplay} with any of the $\dot{x}$-free update laws,~\eqref{UpdateLawICL} or~\eqref{UpdateLawNameWärNett}. Our derivation of such upper bounds, however, necessitates additional smoothness of the regressor:
\begin{assumption}\label{assumption:ExpertKnowledge}
    The required continuity of the regressor $Y$ is to be locally Lipschitz additionally. Furthermore, let $\Xi_\mathrm{t}(B)=\{x_\mathrm{t}(s)+r:s\in[0,T_\mathrm{p}],\, r\in\mathbb{R}^n,\, \|r\|\leq B\}\subset \mathbb{R}^n$ denote the tube of radius $B\geq 0$ around the training segment. Due to $x_\mathrm{t}\in C^1([0,T_\mathrm{p}])$, this tube is compact for each $B\geq 0$. Thereby, $Y$ is Lipschitz on $\Xi_\mathrm{t}(B)$ and there exists $L_\mathrm{t}(B)> 0$ with $\sigma_{\max}(Y(\xi_\mathrm{r})-Y(\xi))\leq L_\mathrm{t}(B)\|\xi-\xi_\mathrm{r}\|$ for all $\xi_\mathrm{r},\xi\in\Xi_\mathrm{t}(B)$.
\end{assumption}

\begin{theorem}\label{theorem:EmployingExpertKnowledge} Consider the closed loop of the plant~\eqref{Plant} and adaptive control~\eqref{AdaptiveControl} under reference replay~\eqref{Replay} where the data-acquisition instants fulfill the replay requirement~\eqref{DataReplay} and the Assumption~\ref{assumption:ExpertKnowledge} applies. For any constant $\bar{V}\geq V(\cdot)|_{t=0}$, the following two statements hold: First, if using~\eqref{UpdateLawICL} as update law while the offline-excitation test $\lambda_{\,\mathrm{t}}^{\!\mathcal{YY}}>0$ is met, then the first occurrence of FE satisfies $T^\star\leq l_{\mathrm{ICL}}^\star T_\mathrm{p}+T_{\max}$~with 
\begin{align}\label{UpperBoundICL}
    l_{\mathrm{ICL}}^\star&=\Big\lfloor\tfrac{\bar{V}T_{\max}}{\lambda_{\min}(K)}\Big(\tfrac{2(\sum_{\tau\in\mathcal{T}_\mathrm{t}}\sigma_{\max}(\mathcal{Y}_\mathrm{t}(\tau)))L_\mathrm{t}({\scriptscriptstyle\sqrt{2\bar{V}}})}{\lambda_{\,\mathrm{t}}^{\!\mathcal{YY}}}\Big)^{\!2}\Big\rfloor,
\end{align} 
i.e.\ FE then happens during or before the $l_{\mathrm{ICL}}^\star\!$th replay~cycle. Second, if $(\ref{UpdateLawNameWärNett})\!$ serves as update law and the offline-excitation test $\lambda_{\,\mathrm{t}}^{\!\mathcal{W}}\!>0$ holds, then the first FE-occurrence happens at $T^\star\leq l^\star T_\mathrm{p}+T_{\max}$ before or during the $l^\star\!$th replay cycle where
\begin{align}\label{UpperBoundNameWärNett}
    l^\star&=\Big\lfloor\tfrac{\bar{V}T_{\max}}{\lambda_{\min}(K)}\Big(\tfrac{(M_\mathrm{t}+M_\mathrm{x}({\scriptscriptstyle\sqrt{2\bar{V}}}))|\mathcal{T}_\mathrm{t}|L_\mathrm{t}({\scriptscriptstyle\sqrt{2\bar{V}}})}{\lambda_{\,\mathrm{t}}^{\!\mathcal{W}}}\Big)^{\!2}\Big\rfloor.
\end{align}
The values $M_\mathrm{t}=\sup_{s\in[0,T_\mathrm{p}]}\{\sigma_{\max}(Y(x_\mathrm{t}(s))\}$ and $M_\mathrm{x}({\scriptstyle\sqrt{2\bar{V}}})=\sup_{\xi\!\in\!\Xi_\mathrm{t}({\scriptscriptstyle\sqrt{2\bar{V}}})}\{\sigma_{\max}(Y(\xi))\}$ in~\eqref{UpperBoundNameWärNett} are finite per the respective compactness of $[0,T_\mathrm{p}]$, $\Xi_\mathrm{t}({\scriptstyle\sqrt{2\bar{V}}})$ since $Y$ is continuous. No-tably, $l_{\mathrm{ICL}}^\star$ and $l^\star$ are arbitrarily small in $\mathbb{N}_0$ by choice of~$K$.
\end{theorem}
\begin{proof} 3a) \textit{Finite contingent of $L^2$-energy during replay.}
Since the time derivative (\eqref{dotV_ICL} for \eqref{UpdateLawICL} and \eqref{dotV_NameWärNett} for \eqref{UpdateLawNameWärNett}) of the Lyapunov function admits $\dot{V}(\cdot)|_{t}\leq -\lambda_{\min}(K)\|e(t)\|^2$ for all $t\in[0,T^\star)$, we obtain $V(\cdot)|_{t<T^\star}-\bar{V}\leq V(\cdot)|_{t<T^\star}-V(\cdot)|_0\leq \int_{0}^{t}\!\!-\lambda_{\min}(K)\|e(s)\|^2\mathrm{d}s$. Thus, $V(\cdot)\geq 0$ reveals the fact that the tracking error is subject to a finite contingent of $L^2$-energy via $\forall t\in[0,T^\star): \int_{0}^{t}\|e(s)\|^2\mathrm{d}s\leq\bar{V}/\lambda_{\min}(K)$, which can (at most) be consumed during the replay phase up until $T^\star$.

3b) \textit{Tube-membership during replay.} We have $x_\mathrm{r}(t)\in\Xi_\mathrm{t}({\scriptstyle\sqrt{ 2\bar{V}}})$ for all $t\in[0,T^\star)$
since $\Xi_\mathrm{t}(\cdot)$ are tubes along the path of $x_\mathrm{r}|_{t<T^\star}$ during replay. Further, for all $t\in[0,T^\star)$ holds $\tfrac{1}{2}\|e(t)\|^2\leq V(\cdot)|_t\leq \bar{V}$ by virtue of $\dot{V}|_t\leq 0$, such that $\|x(t)-x_\mathrm{r}(t)\|\leq {\scriptstyle\sqrt{ 2\bar{V}}}$. Thereby, due to $\Xi_\mathrm{t}({\scriptstyle\sqrt{ 2\bar{V}}})$ containing everything that is inside a ball of radius ${\scriptstyle\sqrt{ 2\bar{V}}}$ around $x_\mathrm{r}(t)|_{t<T^\star}$, also $x(t)\in\Xi_\mathrm{t}({\scriptstyle\sqrt{ 2\bar{V}}})$ for all $t\in[0,T^\star)$. Assumption~\ref{assumption:ExpertKnowledge} thus gives $\forall t\in[0,T^\star):\sigma_{\max}(Y(x_\mathrm{r}(t))-Y(x(t)))\leq L_\mathrm{t}({\scriptstyle\sqrt{ 2\bar{V}}})\|e(t)\|$.

\textit{Proof strategy}---We show that each replay cycle in which FE does not occur consumes a lower-bounded portion of the finite $L^2$-energy contingent that is available for the replay phase, thus limiting the number of FE-absent replay cycles.

3c) \textit{Bounding $T^\star$ under~\eqref{UpdateLawICL} as update law.}  Let $\lambda_{\,\mathrm{t}}^{\!\mathcal{YY}}>0$ be met. We consider any $l\in\mathbb{N}_0$ with $t=T_{\max}+lT_\mathrm{p}$ such that $t<T^\star$. Then, the inequality \eqref{AufhängerTheoremICL} applies. Further, $\lambda^{\!\mathcal{YY}}(t)=0$ due to $t<T^\star$ and thus, we find a $\vartheta\in\mathbb{R}^{p}$ with $\|\vartheta\|=1$ such that the LHS of \eqref{AufhängerTheoremICL} is zero, cf.\ \eqref{lminICL}. That is, we obtain
\begin{align}
    0&\geq{\scriptstyle\lambda_{\,\mathrm{t}}^{\!\mathcal{YY}} - 2 \underset{\tau_\mathrm{t}\in \mathcal{T}_\mathrm{t}}{\sum}\sigma_{\max}(\mathcal{Y}_\mathrm{t}(\tau_\mathrm{t}))\left\|\int_{\tau_\mathrm{t}+lT_\mathrm{p}-\Delta t}^{\tau_\mathrm{t}+lT_\mathrm{p}}\left(Y(x_\mathrm{r}(s))-Y(x(s))\right)\vartheta\mathrm{d}s\right\|}\nonumber\allowdisplaybreaks\\
    &\mkern-9mu\geq {\scriptstyle\lambda_{\,\mathrm{t}}^{\!\mathcal{YY}} - 2 \big(\underset{\tau_\mathrm{t}\in \mathcal{T}_\mathrm{t}}{\sum}\sigma_{\max}(\mathcal{Y}_\mathrm{t}(\tau_\mathrm{t}))\big)\int_{lT_\mathrm{p}}^{lT_\mathrm{p}+T_{\max}}\sigma_{\max}(Y(x_\mathrm{r}(s))-Y(x(s)))\mathrm{d}s}\nonumber\allowdisplaybreaks\\
    &\mkern-9mu\overset{\scriptscriptstyle \text{3b)}}{\geq} {\scriptstyle\lambda_{\,\mathrm{t}}^{\!\mathcal{YY}} - 2 \big(\underset{\tau_\mathrm{t}\in \mathcal{T}_\mathrm{t}}{\sum}\sigma_{\max}(\mathcal{Y}_\mathrm{t}(\tau_\mathrm{t}))\big) L_\mathrm{t}({\scriptscriptstyle\sqrt{2\bar{V}}})\int_{lT_\mathrm{p}}^{lT_\mathrm{p}+T_{\max}}\|e(s)\|\mathrm{d}s} \nonumber\allowdisplaybreaks\\
    \label{AbhängerICL}
    &\mkern-9mu\geq {\scriptstyle\lambda_{\,\mathrm{t}}^{\!\mathcal{YY}} - 2 \big(\underset{\tau_\mathrm{t}\in \mathcal{T}_\mathrm{t}}{\sum}\sigma_{\max}(\mathcal{Y}_\mathrm{t}(\tau_\mathrm{t}))\big) L_\mathrm{t}({\scriptscriptstyle\sqrt{2\bar{V}}})\big(T_{\max}\int_{lT_\mathrm{p}}^{lT_\mathrm{p}+T_{\max}}\|e(s)\|^2\mathrm{d}s\big)^{\!1/2}}\!\!
\end{align}
from~\eqref{AufhängerTheoremICL}. Reorganizing~\eqref{AbhängerICL} yields $\!\int_{lT_\mathrm{p}}^{lT_\mathrm{p}+T_{\max}}\!\!\|e(s)\|^2\mathrm{d}s\!\geq\! b^\star_\mathrm{ICL}$ with\footnote{$\lambda_{\,\mathrm{t}}^{\!\mathcal{YY}}>0$, via~\eqref{TrainingICL}, guarantees the denominator of $b^\star_\mathrm{ICL}$ to be positive.} $\!b^\star_\mathrm{ICL}\!=\!(\lambda_{\,\mathrm{t}}^{\!\mathcal{YY}})^2(2(\sum_{\tau_\mathrm{t}\in\mathcal{T}_\mathrm{t}}\sigma_{\max}(\mathcal{Y}_\mathrm{t}(\tau_\mathrm{t})))L_\mathrm{t}({\scriptstyle\sqrt{2\bar{V}}}))^{-2}T_{\max}^{-1}\!>0$. 
Therefore, any $l'\in\mathbb{N}_0$ with $l'T_\mathrm{p}+T_{\max}<T^\star$ must admit $\!\int_{0}^{l'T_\mathrm{p}+T_{\max}}\!\|e(s)\|^2\mathrm{d}s\geq \sum_{l=0}^{l'}\int_{lT_\mathrm{p}}^{lT_\mathrm{p}+T_{\max}}\|e(s)\|^2\mathrm{d}s\geq (l'+1)b^\star_\mathrm{ICL}$. In view of 3a), this makes $\bar{V}/\lambda_{\min}(K)\geq (l'+1)b^\star_\mathrm{ICL}$ necessary for $l'T_\mathrm{p}+T_{\max}<T^\star$ when $l'\in\mathbb{N}_0$. Accordingly, all $l'\in\mathbb{N}_0$ with $l'>(\bar{V}/\lambda_{\min}(K))(b_\mathrm{ICL}^\star)^{-1}-1$ admit $l'T_\mathrm{p}+T_{\max}\geq T^{\star}$, the smallest of which is $l_{\mathrm{ICL}}^\star$ from~\eqref{UpperBoundICL}, as claimed.

3d) \textit{Upper bound on $T^\star$ under~\eqref{UpdateLawNameWärNett} serving as update law.} Suppose, $\lambda_{\,\mathrm{t}}^{\!\mathcal{W}}\!>0$. For any $l\in\mathbb{N}_0$ with $t=T_{\max}+lT_\mathrm{p}$ such that $t<T^\star$, the following holds: The inequality \eqref{AufhängerTheoremNameWärNett} applies and, due to $\lambda^{\!\mathcal{W}}(t)=0$ as per $t<T^\star$, we can find a $\vartheta\in\mathbb{R}^p$ with $\|\vartheta\|=1$ such that its LHS is zero. For this $\vartheta$, the terms on the RHS of \eqref{AufhängerTheoremNameWärNett} admit $\|(Y(x_\mathrm{r}(s))+Y(x(s)))\vartheta\|\leq M_\mathrm{t}+M_\mathrm{x}({\scriptstyle\sqrt{2\bar{V}}})$ for all $s\in [\tau_\mathrm{t}+lT_\mathrm{p}-\Delta t,\tau_\mathrm{t}+lT_\mathrm{p}]$, every $\tau_\mathrm{t}\in\mathcal{T}_\mathrm{t}$ (cf. 3b)). Thus, we may lower bound the RHS of \eqref{AufhängerTheoremNameWärNett} and~obtain
\begin{align}
    0&\geq {\scriptstyle\lambda_{\,\mathrm{t}}^{\!\mathcal{W}}-\underset{\tau_\mathrm{t}\in\mathcal{T}_\mathrm{t}}{\sum}\int_{\tau_\mathrm{t}+lT_\mathrm{p}-\Delta t}^{\tau_\mathrm{t}+lT_\mathrm{p}}(M_\mathrm{t}+M_\mathrm{x}({\scriptscriptstyle\sqrt{2\bar{V}}}))\|\left(Y(x_\mathrm{r}(s))-Y(x(s))\right)\vartheta\| \mathrm{d}s }\allowdisplaybreaks\nonumber\\
    &\geq {\scriptstyle\lambda_{\,\mathrm{t}}^{\!\mathcal{W}}-(M_\mathrm{t}+M_\mathrm{x}({\scriptscriptstyle\sqrt{2\bar{V}}}))|\mathcal{T}_\mathrm{t}|\int_{\tau_\mathrm{t}+lT_\mathrm{p}-\Delta t}^{\tau_\mathrm{t}+lT_\mathrm{p}}\sigma_{\max}(\left(Y(x_\mathrm{r}(s))-Y(x(s))\right)) \mathrm{d}s }\allowdisplaybreaks\nonumber\\
    &\overset{\scriptscriptstyle \text{3b)}}{\geq} {\scriptstyle\lambda_{\,\mathrm{t}}^{\!\mathcal{W}}-(M_\mathrm{t}+M_\mathrm{x}({\scriptscriptstyle\sqrt{2\bar{V}}}))|\mathcal{T}_\mathrm{t}|L_\mathrm{t}({\scriptscriptstyle\sqrt{2\bar{V}}})\int_{lT_\mathrm{p}}^{lT_\mathrm{p}+T_{\max}}\|e(s)\| \mathrm{d}s }\allowdisplaybreaks\nonumber
    \\
    &\geq {\scriptstyle\lambda_{\,\mathrm{t}}^{\!\mathcal{W}}-(M_\mathrm{t}+M_\mathrm{x}({\scriptscriptstyle\sqrt{2\bar{V}}}))|\mathcal{T}_\mathrm{t}|L_\mathrm{t}({\scriptscriptstyle\sqrt{2\bar{V}}})\big(T_{\max}\int_{lT_\mathrm{p}}^{lT_\mathrm{p}+T_{\max}}\|e(s)\|^2 \mathrm{d}s\big)^{\!1/2}}.
\end{align}
For any such $l\in\mathbb{N}_0$, this reveals $\int_{lT_\mathrm{p}}^{lT_\mathrm{p}+T_{\max}}\!\!\|e(s)\|^2\mathrm{d}s\!\geq\! b^\star$ with\footnote{Due to $\lambda_{\,\mathrm{t}}^{\!\mathcal{W}}\!>0$, the denominator of $b^\star$ must be positive, cf.\ \eqref{TrainingNameWärNett}.} $b^\star=(\lambda_{\,\mathrm{t}}^{\!\mathcal{W}})^2((M_\mathrm{t}+M_\mathrm{x}({\scriptstyle\sqrt{2\bar{V}}}))|\mathcal{T}_\mathrm{t}|L_\mathrm{t}({\scriptstyle\sqrt{2\bar{V}}}))^{-2}T_{\max}^{-1}\!>0$. Thus, we may proceed analogously to 3c) by noting that any $l''\in\mathbb{N}_0$ with $l''T_\mathrm{p}+T_\mathrm{\max}<T^\star$ satisfies $\int_{0}^{l''T_\mathrm{p}+T_{\max}}\|e(s)\|^2\mathrm{d}s\geq \sum_{l=0}^{l''}\int_{lT_\mathrm{p}}^{lT_\mathrm{p}+T_{\max}}\!\!\|e(s)\|^2\mathrm{d}s\!\geq\! (l''\!\!+\!1)b^\star\!$. In view of 3a), this~shows that $\bar{V}/\lambda_{\min}(K)\!\geq\! (l''\!\!+\!1)b^\star$ is necessary for any $l''\!\in\mathbb{N}_0$ to admit $l''T_\mathrm{p}\!+\!T_{\max}<T^\star\!$. Accordingly, any $l''\!\in\mathbb{N}_0$ that~violates this necessity via $l''\!>(\bar{V}/\lambda_{\min}(K))(b^\star)^{-1}\!-\!1$ satisfies $l''T_\mathrm{p}+T_{\max}\geq T^\star\!$, where $l^\star\!$ from~\eqref{UpperBoundNameWärNett} is the smallest of these $l''\!$.
\end{proof}

\textit{Leveraging expert knowledge---}In case knowledge is available about the region in which the unknown parameters $\theta$ are located so that an upper bound $\bar{V}$ on the unknown initial value $V(\cdot)|_{t=0}$ of the Lyapunov function from \eqref{LyapFcn} can be found, then the values $l_{\mathrm{ICL}}^\star$, $l^\star$ from~\eqref{UpperBoundICL},~\eqref{UpperBoundNameWärNett} are computable \textit{offline}; thus indicating \textit{a priori} the number of replay cycles in which FE must occur under the presumptions of Theorem~\ref{theorem:EmployingExpertKnowledge}. Interestingly, then pursuing a high-gain strategy for~\eqref{AdaptiveControl} with
\begin{align}
\lambda_{\min}(K)>  \left\{\begin{matrix}
        {\scriptstyle \bar{V}T_{\max}}\Big(\tfrac{2(\sum_{\tau\in\mathcal{T}_\mathrm{t}}\sigma_{\max}(\mathcal{Y}_\mathrm{t}(\tau)))L_\mathrm{t}({\scriptscriptstyle\sqrt{2\bar{V}}})}{\lambda_{\,\mathrm{t}}^{\!\mathcal{YY}}}\Big)^{\!2}:&\!\!\!\!{\scriptstyle \text{\parbox{1.6cm}{%
        \linespread{0.5}\selectfont
        $\lambda_{\,\mathrm{t}}^{\!\mathcal{YY}}\!>\!0$ and $\dot{\hat{\theta}}$ as per \eqref{UpdateLawICL}}}}\\
        {\scriptstyle\bar{V}T_{\max}}\Big(\tfrac{(M_\mathrm{t}+M_\mathrm{x}({\scriptscriptstyle\sqrt{2\bar{V}}}))|\mathcal{T}_\mathrm{t}|L_\mathrm{t}({\scriptscriptstyle\sqrt{2\bar{V}}})}{\lambda_{\,\mathrm{t}}^{\!\mathcal{W}}}\Big)^{\!2}:&\!\!\!\!{\scriptstyle \text{\parbox{1.6cm}{%
        \linespread{0.5}\selectfont
        $\lambda_{\,\mathrm{t}}^{\!\mathcal{W}}\!>\!0$ and $\dot{\hat{\theta}}$ as per \eqref{UpdateLawNameWärNett}}}}
    \end{matrix}\right.,
\end{align}
under fulfillment of the respective offline excitation test and validity of Assumption~\ref{assumption:ExpertKnowledge}, enforces FE during the $0$th replay cycle at $T^\star\in[0, T_{\max}]$, cf.\ Theorem \ref{theorem:EmployingExpertKnowledge}. That is, FE~occurrence\linebreak then takes no longer than it takes during runtime to experience the full excitation of the training segment once: $T_{\max}$.

\subsection{$\!$Comparative view on offline-excitation-test data efficiency}
The a-priori-verifiable conditions needed in Theorems \ref{Theorem:ICL}, \ref{Theorem:NameWärNett}, \ref{theorem:EmployingExpertKnowledge} to guarantee finite-time FE-injection under \eqref{Replay}, \eqref{DataReplay} ask for positive signs of the offline-computable, minimal eigenvalues $\lambda_{\,\mathrm{t}}^{\!\mathcal{YY}}$, $\lambda_{\,\mathrm{t}}^{\!\mathcal{W}}$. Notably, success or failure of these offline excitation tests for a jointly given training segment are connected to each other; as is evident by the following~result:
\begin{theorem}\label{theorem:DataEfficiency}
   It holds $\sum_{\tau_\mathrm{t}\in\mathcal{T}_\mathrm{t}}\mathcal{Y}_\mathrm{t}(\tau_\mathrm{t})^\top\mathcal{Y}_\mathrm{t}(\tau_\mathrm{t})\preceq \Delta t \sum_{\tau_\mathrm{t}\in\mathcal{T}_\mathrm{t}} \mathcal{W}_\mathrm{t}(\tau_\mathrm{t})$.
\end{theorem}
\begin{proof}
    For any $\vartheta\in\mathbb{R}^p$, define $h_\vartheta(\cdot)=Y(x_\mathrm{t}(\cdot))\vartheta\in\mathbb{R}^n$ and denote its rows as $h_\vartheta^i(\cdot)\in\mathbb{R}$, $i\in\{1,\ldots,n\}$, i.e.\ $h_\vartheta(\cdot)=(\begin{smallmatrix}
        h_\vartheta^1(\cdot)&\ldots&h_\vartheta^n(\cdot)
    \end{smallmatrix})\!^\top$. Observe that $\|\int_{\tau_\mathrm{t}-\Delta t}^{\tau_\mathrm{t}} h_\vartheta(s)\mathrm{d}s\|^2=\sum_{i=1}^n(\int_{\tau_\mathrm{t}-\Delta t}^{\tau_\mathrm{t}} h_\vartheta^i(s)\mathrm{d}s)^2\leq \sum_{i=1}^n\int_{\tau_\mathrm{t}-\Delta t}^{\tau_\mathrm{t}} 1^2\mathrm{d}s \int_{\tau_\mathrm{t}-\Delta t}^{\tau_\mathrm{t}} h_\vartheta^i(s)^2\mathrm{d}s= \Delta t \int_{\tau_\mathrm{t}-\Delta t}^{\tau_\mathrm{t}} \sum_{i=1}^n h_\vartheta^i(s)^2\mathrm{d}s=\Delta t \int_{\tau_\mathrm{t}-\Delta t}^{\tau_\mathrm{t}} \|h_\vartheta(s)\|^2\mathrm{d}s$ for any $\tau_\mathrm{t}\in\mathcal{T}_\mathrm{t}$, by Cauchy–Schwarz' inequality justifying the less-or-equal sign therein. Thereby, we obtain $\vartheta^\top \mathcal{Y}_\mathrm{t}(\tau_\mathrm{t})^\top \mathcal{Y}_\mathrm{t}(\tau_\mathrm{t})\vartheta=\|\int_{\tau_\mathrm{t}-\Delta t}^{\tau_\mathrm{t}} h_\vartheta(s)\mathrm{d}s\|^2\leq \Delta t \int_{\tau_\mathrm{t}-\Delta t}^{\tau_\mathrm{t}} \|h_\vartheta(s)\|^2\mathrm{d}s=\Delta t \vartheta^\top \mathcal{W}_\mathrm{t}(\tau_\mathrm{t})\vartheta$, which yields the claim via $\mathcal{Y}_\mathrm{t}(\tau_\mathrm{t})^\top\mathcal{Y}_\mathrm{t}(\tau_\mathrm{t})\preceq \Delta t \mathcal{W}_\mathrm{t}(\tau_\mathrm{t})$.
\end{proof}

In view of~\eqref{TrainingICL},~\eqref{TrainingNameWärNett}, Theorem~\ref{theorem:DataEfficiency} thus leads to
\begin{align}\label{PathsAreExciting}
    \lambda_{\,\mathrm{t}}^{\!\mathcal{YY}}>0&\Longrightarrow \lambda_{\,\mathrm{t}}^{\!\mathcal{W}}>0,\\
    \lambda_{\,\mathrm{t}}^{\!\mathcal{W}}=0&\Longrightarrow \lambda_{\,\mathrm{t}}^{\!\mathcal{YY}}=0.
\end{align}
That is, whenever the offline excitation test is fulfilled w.r.t.\ the data-augmented update law~\eqref{UpdateLawICL} of ICL, then it also~holds w.r.t.\ the data-augmentation of the update law \eqref{UpdateLawNameWärNett}. Conversely, ICL cannot pass its offline excitation test for any training segment on which the data-acquisition method of \eqref{UpdateLawNameWärNett} fails to collect exciting data. Thus, the (in~Section~\ref{section:PRELIMINARIES} highlighted) computational complexity of the data-acquisi-\linebreak tion method \eqref{JaAnscheinendIstDasSo}-\eqref{HachJa} to facilitate the data-augmentation of \eqref{UpdateLawNameWärNett} achieves being more efficient than ICL's data acquisition \eqref{FTCTrickery}, \eqref{Rattenschwanz}. Therefore, it is preferable to use~\eqref{UpdateLawNameWärNett} as update law compared to~\eqref{UpdateLawICL} whenever the regressor fulfills Assumption~\ref{GradientAssumption}. This originates from the fact that under Assumption~\ref{GradientAssumption}, data acquisition for~\eqref{UpdateLawNameWärNett} is possible via \eqref{PathInGradientFields}, \eqref{ForDiscussionsSake}, \eqref{HachJa}, which\newpage\noindent is equally as demanding as \eqref{FTCTrickery}, \eqref{Rattenschwanz} for~\eqref{UpdateLawICL} while retaining superiority in data efficiency. In case the regressor does not fulfill Assumption~\ref{GradientAssumption}, then the choice of update law can be made in dependency on the specific training segment: If ICL passes its offline excitation test $ \lambda_{\,\mathrm{t}}^{\!\mathcal{YY}}>0$, then one may use~\eqref{UpdateLawICL} as update law in order to enjoy computationally light data acquisition via~\eqref{FTCTrickery}, \eqref{Rattenschwanz}. Should ICL's offline excitation test fail while $\lambda_{\,\mathrm{t}}^{\!\mathcal{W}}>0$ holds then one may resort to the update law~\eqref{UpdateLawNameWärNett} while accepting computationally demanding data acquisition with \eqref{JaAnscheinendIstDasSo}-\eqref{HachJa}. Otherwise, another training segment must be tested offline for its excitation.

\subsection{Contextualization of reference replay w.r.t.\ \cite{gallegos_relaxed_2024}: From a-priori FE-existence guarantees to upper time-to-FE bounds}\label{subsection:RRvsTE}
To the best of the authors' knowledge, a-priori guarantees for closed-loop FE-occurrence appear in the prior literature only in (\cite{gallegos_relaxed_2024}, Thm.\ 5.6). Therein, the integral-memory estimator (\cite{gallegos_relaxed_2024}, Eq.\ 5.12) is shown to collect FE-data in finite, unknown time $T^\star\!<\infty$, provided the desired reference~tra-jectory is replaced by a signal $\eta$ on which the regressor is PE until FE occurs. The reference-replay mechanism~\eqref{Replay} works similarly since it replaces the desired reference trajectory~$\!x_{\mathrm{d},T^\star}\!\!$ with the $T_\mathrm{p}$-periodic signal $x_\mathrm{t}(\cdot-\lfloor\tfrac{\cdot}{T_\mathrm{p}}\rfloor T_\mathrm{p}\rfloor)$ until FE occurs. Notably, the regressor $Y$ evaluated on this signal is rendered PE whenever any of the two offline excitation tests holds,~see:
\begin{definition} (\cite{gallegos_relaxed_2024}, Def.\ 3.1) A signal $m:[0,\infty)\!\to\!\mathbb{R}^{n\times p}\!$ is PE if $\epsilon'\!,\delta'\!>0$ exist such that $\forall t\geq 0: \int_{t}^{t+\epsilon'}\!m(s)^\top m(s)\mathrm{d}s\geq \delta' \mathrm{I}_p$.
\end{definition}
 \begin{corollary}\label{corollary:PEGenerator}
Consider $\eta$ with $\eta(t)=x_\mathrm{t}(t-\lfloor t/T_\mathrm{p}\rfloor T_\mathrm{p})$ for all $t\! \geq\! 0$. If the training segment is such that one of the offline~ex-citation tests  holds, i.e.\ if $\lambda_{\,\mathrm{t}}^{\!\mathcal{YY}}\!\!>\!0$ or $\!\lambda_{\,\mathrm{t}}^{\!\mathcal{W}}\!\!>\!0$, then $Y\!(\eta)$ is~PE.  
\end{corollary}
\begin{proof}
    Let $\lambda_{\,\mathrm{t}}^{\!\mathcal{YY}}\!\!>\!0$ or $\lambda_{\,\mathrm{t}}^{\!\mathcal{W}}\!\!>\!0$. Due to \eqref{PathsAreExciting}, we thus have\linebreak  $\sum_{\tau_\mathrm{t}\in\mathcal{T}_\mathrm{t}}\int_{I(\tau_\mathrm{t})}D_\mathrm{t}(s)\mathrm{d}s\succ 0$, wherein $D_\mathrm{t}(s)=Y(x_\mathrm{t}(s))^\top Y(x_\mathrm{t}(s))$ and $I(\tau_\mathrm{t})=[\tau_\mathrm{t}-\Delta t,\tau_\mathrm{t}]$ serve as shortcuts for brevity. Now, if one assumes $\int_{\cup_{\tau_\mathrm{t}\in\mathcal{T}_\mathrm{t}}I(\tau_{\mathrm{t}})}\!D_\mathrm{t}(s)\mathrm{d}s\!\nsucc\! 0$, then there must exist a\linebreak $0\!\neq\!\vartheta\!\in\!\mathbb{R}^p\!$ such that $\int_{\cup_{\tau_\mathrm{t}\in\mathcal{T}_\mathrm{t}}\!I(\tau_{\mathrm{t}})}\vartheta^\top D_\mathrm{t}(s)\vartheta\mathrm{d}s\!=\!0$. Because the\linebreak  dyadic products $\vartheta^\top D_\mathrm{t}(\cdot)\vartheta$ therein are continuous and~non-negative, this leads to $(\vartheta^\top D_\mathrm{t}(\cdot)\vartheta)|_{I(\tau_\mathrm{t})}\equiv 0$ for all $\tau_\mathrm{t}\in\mathcal{T}_\mathrm{t}$.\linebreak That is, non-positive-definiteness would produce~a~contradic-\linebreak tion, $\forall\tau_\mathrm{t}\!\in\!\mathcal{T}_\mathrm{t}:\vartheta^\top \int_{I(\tau_\mathrm{t})}\!D_\mathrm{t}(s)\mathrm{d}s\vartheta\!=\!0$, to $\sum_{\tau_\mathrm{t}\in\mathcal{T}_\mathrm{t}}\int_{I(\tau_\mathrm{t})}D_\mathrm{t}(s)\mathrm{d}s\succ 0$. This in turn proves $\int_{\cup_{\tau_\mathrm{t}\in\mathcal{T}_\mathrm{t}}I(\tau_{\mathrm{t}})}\!D_\mathrm{t}(s)\mathrm{d}s$ to be positive~definite. Thus, $\int_{[0,T_\mathrm{p}]}D_\mathrm{t}(s)\mathrm{d}s\succeq \int_{\cup_{\tau_\mathrm{t}\in\mathcal{T}_\mathrm{t}}I(\tau_{\mathrm{t}})}D_\mathrm{t}(s)\mathrm{d}s\succ 0$ must hold since $ [0,T_\mathrm{p}]\supseteq \cup_{\tau_\mathrm{t}\in\mathcal{T}_\mathrm{t}}I(\tau_{\mathrm{t}})$ is guaranteed via $[\Delta t, T_\mathrm{p}]\supset \mathcal{T}_\mathrm{t}$. Accordingly, $\delta'=\lambda_{\min}(\int_{[0,T_\mathrm{p}]}D_\mathrm{t}(s)\mathrm{d}s)$ is positive. Let $D_\mathrm{\eta}(\cdot)=Y(\eta(\cdot))^\top Y(\eta(\cdot))$ and $\epsilon'=T_\mathrm{p}$. Thus, $\epsilon'$ is also positive.~Note that $\forall t\!\geq\! 0: \int_{t}^{t+\epsilon'}\!D_\mathrm{\eta}(s)\mathrm{d}s=\int_{[0,T_\mathrm{p}]}\!D_\mathrm{t}(s)\mathrm{d}s$ since $D_\mathrm{\eta}$ is the~$T_\mathrm{p}$-periodic repetition of $D_\mathrm{t}$. Thereby, $\epsilon',\delta'>0$ fit the claimed PE requirement with $\forall t\geq 0:\lambda_{\min}(\int_{t}^{t+\epsilon'}D_\mathrm{\eta}(s)\mathrm{d}s)=\delta'$.
\end{proof} 

Thereby, the reference-replay mechanism~\eqref{Replay} can be understood as a generator of reference trajectories for which the estimator from~\cite{gallegos_relaxed_2024} is proven therein to collect FE-data in finite, unknown time whenever any of the offline excitation tests, $\lambda_{\,\mathrm{t}}^{\!\mathcal{YY}}\!>0$ or $\lambda_{\,\mathrm{t}}^{\!\mathcal{W}}\!>0$, holds for a given training segment.

Despite this similarity in the generation of reference~trajec-tories, the findings herein on a-priori FE-verifiability differ from those of \cite{gallegos_relaxed_2024} in two structural aspects: First, the~replay-\newpage\noindent requirement \eqref{DataReplay} synchronizes the data acquisition with the replay cycle. In Theorems \ref{Theorem:ICL}, \ref{Theorem:NameWärNett} this allows us to establish the central existence result $T^\star<\infty$ of (\cite{gallegos_relaxed_2024}, Thm.\ 5.6) also for CL where data is acquired at discrete time instances~$\mathcal{T}(\cdot)$, rather than by incorporating complete trajectories of closed-loop signals during $[0,T^\star)$ into the estimation, as is done in \cite{gallegos_relaxed_2024}. Second, the proof-machinery of Theorems \ref{Theorem:ICL}, \ref{Theorem:NameWärNett} enables us in Theorem~\ref{theorem:EmployingExpertKnowledge} to quantify the preceding existence results by providing upper bounds on $T^\star$---which is absent in~\cite{gallegos_relaxed_2024}.

\vspace{-.5mm}\section{Simulation}\label{section:SIMULATION}\vspace{-.5mm}
\begin{figure}
\vspace{-0mm}
\centerline{\includegraphics[scale=0.59]{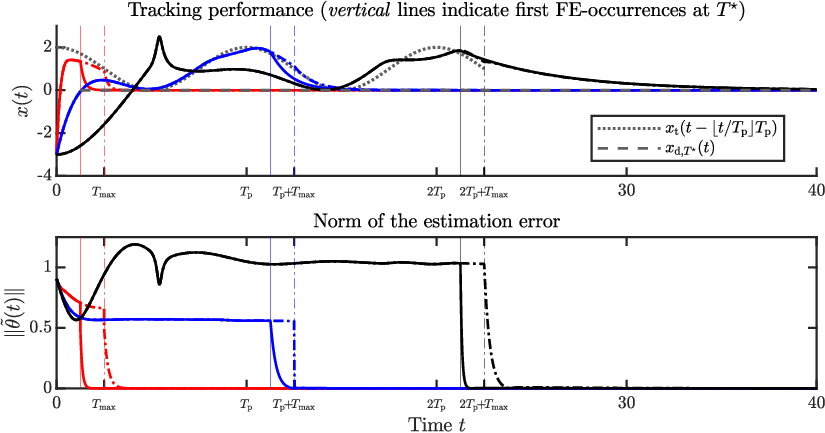}}
\vspace{-1.75mm}
\caption{Closed-loop performance under reference replay when increasing the control gain (\textit{black} $K\!=\!.2$, \textit{blue} $K\!=\!1$, \textit{red} $K\!=\!5$). \textit{Dash-dotted} lines use \eqref{UpdateLawICL} as update law; \textit{solid} lines use \eqref{UpdateLawNameWärNett}. FE-data is collected during all experiment-runs, subsequent to which tracking and estimation errors~con-verge toward zero. Notably, higher control gains cause~FE to occur earlier, down to the 0th replay cycle. Moreover, at any gain, the update law~\eqref{UpdateLawNameWärNett} collects FE-data earlier than~\eqref{UpdateLawICL}.}
\label{fig:Replay}
\vspace{-4mm}
\end{figure}

Fig.\ \ref{fig:Replay} shows the simulation results of a plant with $p=2$\linebreak parameters $\theta=-({\scriptstyle \tfrac{1}{10},\ \tfrac{9}{10}})^\top$ (unknown to the controller) and $1$-dimensional state $x(t)\in\mathbb{R}$, $x(0)=-3$ in closed loop with adaptive control~\eqref{AdaptiveControl} and any of the update laws \eqref{UpdateLawICL}, \eqref{UpdateLawNameWärNett} under reference replay~\eqref{Replay} when the control gain $K>0$ is~increased. The regressor is $Y(x)=\Big\{\begin{smallmatrix}
    (1,\ 1):&x<1\\
    (x^2,\ x^3):&x\geq 1
\end{smallmatrix}$ and fulfills Assumption~\ref{GradientAssumption} with $({\scriptstyle \Phi_1(x),\ \Phi_2(x)})=\Big\{\begin{smallmatrix}
    (x-2/3,\ x-3/4):&x<1\\
    (x^3/3,\ x^4/4):&x\geq 1
\end{smallmatrix}$, which we employ to acquire data for~\eqref{UpdateLawNameWärNett} via \eqref{PathInGradientFields}. Further, \eqref{Replay} replays the training segment $x_\mathrm{t}(t)=\cos(\tfrac{2\pi}{T_\mathrm{p}}t)+1$ of length $T_\mathrm{p}=10$ until FE occurs and tracking of the desired reference trajectory $x_{\mathrm{d},T^\star}\equiv 0$ is permitted. Training data is acquired at $\mathcal{T}_\mathrm{t}=\{{\scriptstyle 1.25,\ 2.5}\}$ (i.e., $T_{\max}=2.5$) with $\Delta t=1.25$ serving as buffer. Thus, both offline excitation tests are fulfilled with $\lambda_{\,\mathrm{t}}^{\!\mathcal{YY}}\!=.0366$, $\lambda_{\,\mathrm{t}}^{\!\mathcal{W}}\!=.0922$.
The update laws are parameterized with $\Gamma=\tfrac{1}{10}\mathrm{I}_2$, $k(t)=\Big\{\begin{smallmatrix}
    0:&t\in[0,T^\star)\\
    5000:&t\in [T^\star,\infty)
\end{smallmatrix}\!$~and initialized at $\hat{\theta}(0)=({\scriptstyle 0,\ 0})^{\top}$. Lastly, data is acquired online at $\mathcal{T}(\cdot)$ as per \eqref{ExplicitDataReplay} in synchronization with the replay~cycle. 

First, the fact that FE-data is collected during each of the shown experiments, i.e.\ $T^\star<\infty$, verifies Theorems \ref{Theorem:ICL}, \ref{Theorem:NameWärNett}. Second, Theorem~\ref{theorem:EmployingExpertKnowledge} proves FE to happen within $l^\star_{\mathrm{ICL}},l^\star\to 0$ replay cycles when $K\to\infty$ for \eqref{UpdateLawICL}, \eqref{UpdateLawNameWärNett}, respectively. This phenomenon is also evident in the simulation as Fig.~\ref{fig:Replay} reports FE-occurrence within the $0$th replay cycle $[0,T_\mathrm{p}]$, specifically during $[0,T_{\max}]$, when using a high gain of $K=5$.

\vspace{-.5mm}\section*{Appendix}\vspace{-.5mm}

\begin{lemma}\label{ReplayBisDerArztKommt}
    If $x_\mathrm{r}(t)=x_\mathrm{t}(t-\lfloor t/T_\mathrm{p}\rfloor T_\mathrm{p})$ for all $t\geq 0$ and{\parfillskip=0pt\par}\newpage\noindent $\lim_{t\to\infty}e(t)=0$, then continuity of $Y$ and $x_\mathrm{t}\in C^1([0,T_\mathrm{p}])$~yield
    \begin{align}
        \forall \epsilon>0:\,\exists T(\epsilon)\geq \Delta t:\, \mathcal{S}_{\Delta t}(t)\leq \epsilon\,\forall t\geq T(\epsilon)
    \end{align}
    for the supremum $\mathcal{S}_{\Delta t}$ from~\eqref{SupremeMaxSing}. Therein, $T(\epsilon)\geq T_\mathrm{p}$ w.l.o.g.
\end{lemma}
\begin{proof}
    Let $\Xi_\mathrm{t}(B)=\{x_\mathrm{t}(s)+r:s\in[0,T_\mathrm{p}],\, r\in\mathbb{R}^n,\, \|r\|\leq B\}\subset \mathbb{R}^n$ be the tube of some fixed radius $B>0$ around the path of $x_\mathrm{r}$. Due to the continuity and hence, the boundedness of $x_\mathrm{t}$ on $[0,T_\mathrm{p}]$, this tube is compact and thus, the tube-restricted regressor $Y|_{\Xi_\mathrm{t}(B)}$ is uniformly continuous (UC). We define $S(\xi_\mathrm{r},\xi)=\sigma_{\max}\{Y(\xi_\mathrm{r})-Y(\xi)\}$ for $\xi_\mathrm{r},\xi\in\mathbb{R}^n$. UC of $\Xi_\mathrm{t}(B)$ then guarantees, for any $\epsilon>0$, that $S(\xi_\mathrm{r},\xi)\leq \epsilon$ whenever $\xi_\mathrm{r},\xi\in \Xi_\mathrm{t}(B)$ with $\|\xi_\mathrm{r}-\xi\|\leq \delta_B(\epsilon)$ for some suitable $\delta_B(\epsilon)>0$. As $e(t)\to 0$ for $t\to \infty$, we know that there must exist $T_\mathrm{e}(\min\{\delta_B(\epsilon),B\})\geq 0$ such that, for any $t\geq T_\mathrm{e}(\min\{\delta_B(\epsilon),B\})$, the following holds: $\|x_\mathrm{r}(t)-x(t)\|\leq \min\{\delta_B(\epsilon),B\}$ and thereby, due to the design of $\Xi_\mathrm{t}(B)$ to contain everything that is inside of a ball of radius $B$ around $x_\mathrm{r}(t)\in\Xi_\mathrm{t}(B)$, we also have $x(t)\in\Xi_\mathrm{t}(B)$. That is, we obtain $S(x_\mathrm{r}(s),x(s))\leq \epsilon$ for all $s\geq T_\mathrm{e}(\min\{\delta_B(\epsilon),B\})$. Accordingly, $\mathcal{S}_{\Delta t}(t)\leq \epsilon$ for all $t\geq T_\mathrm{e}(\min\{\delta_B(\epsilon),B\})\!+\!\Delta t=T(\epsilon)\geq \Delta t$.
\end{proof}


\addtolength{\textheight}{-12cm}   

\bibliographystyle{IEEEtran}
\bibliography{Bibliothek}

\end{document}